\documentclass{amsart}
\usepackage{amsfonts}
\usepackage[all]{xy}
\usepackage{amssymb, amscd}

\usepackage{amssymb}
\usepackage{graphicx}
\usepackage{epstopdf}

\usepackage[dvips,bookmarks=false]{hyperref}

\newtheorem{theorem}{Theorem}
\newtheorem{lemma}{Lemma}
\newtheorem{proposition}{Proposition}
\newtheorem{corollary}{Corollary}

\newcommand{\R}{{\Bbb R}}
\newcommand{\C}{{\Bbb C}}
\renewcommand{\H}{{\Bbb H}}
\renewcommand{\O}{{\Bbb O}}
\newcommand{\Z}{{\Bbb Z}}

\renewcommand{\t}{\mathfrak{t}}
\newcommand{\h}{\mathfrak{h}}
\renewcommand{\l}{\mathfrak{l}}
\newcommand{\s}{\mathfrak{s}}
\renewcommand{\k}{\mathfrak{K}}

\newcommand{\g}{\mathfrak{g}}
\newcommand{\p}{\mathfrak{p}}

\newcommand{\m}{\mathfrak{m}}
\newcommand{\so}{\mathfrak{so}}
\renewcommand{\sl}{\mathfrak{sl}}

\newcommand{\End}{{\rm End}}
\newcommand{\GL}{{\rm GL}}

\newcommand{\Aut}{{\rm Aut}}

\newcommand{\SO}{{\rm SO}}
\newcommand{\SL}{{\rm SL}}

\newcommand{\SU}{{\rm SU}}
\newcommand{\x}{{\bf x}}
\newcommand{\y}{{\bf y}}
\renewcommand{\P}{{\bf P}}
\newcommand{\bp}{{\bf p}}

\renewcommand{\c}{ {\bf c}}

\renewcommand{\L}{{\mathcal L}}

\newcommand{\bs}{\bigskip}

\newcommand{\mn}{\medskip\noindent}
\newcommand{\n}{\noindent}
\newcommand{\bn}{\bs\n}

\newcommand{\X}{\times}
\renewcommand{\o}{{\omega}}

\newcommand{\pf}{\n\textsc{Proof. }}
\newcommand{\e}{ {\bf e}}
\renewcommand{\a}{ {\bf a}}
\renewcommand{\b}{ {\bf b}}
\newcommand{\tQ}{\widetilde {Q}}
\newcommand{\tD}{\widetilde {D}}
\newcommand{\tC}{\widetilde {C}}

\title{$G_2$ and the ``rolling distribution''}

\author{ Gil Bor and Richard Montgomery}

\begin{document}
\maketitle

\tableofcontents

\section*{0 Introduction}

\begin{figure}[h]
\centerline{\includegraphics[width=8cm, angle=90]{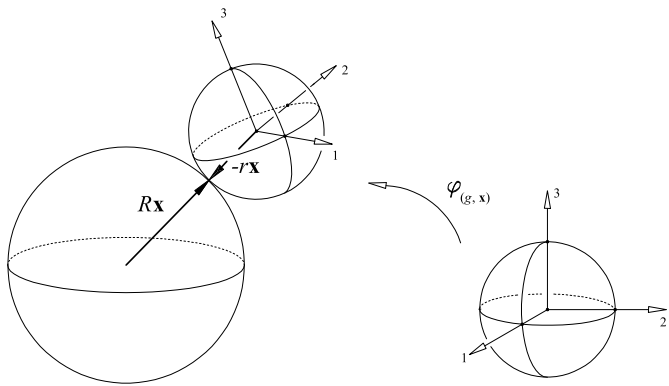}} 
\caption{Rolling a ball on another ball.}\label{roll}
\end{figure}

Consider two balls of different sizes,
rolling on each other, without slipping or spinning. 
The
configuration space for this system is a 5-dimensional manifold
 $Q \cong \SO_3\times
S^2$ on which  the no-slip/no-spin condition defines a rank 2
distribution $D\subset TQ$, the ``rolling-distribution''.

Now $D$  is a non-integrable distribution (unless the balls are of equal size) which has an ``obvious'' 6-dimensional transitive symmetry group 
$\SO_3\times\SO_3$ arising from the isometry groups of each ball, but for balls whose radii are in the ratio 3:1, and only for this ratio, 
something strange happens: the {\it local } symmetry group 
of  the distribution  increases  from  $\SO_3\times\SO_3$  to $G_2$,  a
14-dimensional Lie group. 

More  precisely, let $\g_2$ be the  real  ``split
form"   of the complex $14$-dimensional exceptional simple Lie algebra $\g_{2 \C}$.
There are precisely two connected Lie groups whose Lie algebras
are $\g_2$.  (See Appendix  \ref{twoG2s}.) We choose the one corresponding to the adjoint
representation, and call it $G_2$.  (The other one is $\tilde G_2$,
the universal cover of the one we chose.)
$G_2$ is  a subgroup of $\SO(3,4)$ and its maximal compact subgroup  $K\subset G_2$ is isomorphic to $\SU_2\times \SU_2/\{\pm(1,1)\}$ which   double-covers  $\SO_3\times\SO_3$. 
Let $\tQ=S^3\times S^2$ be 
the universal  cover of $Q$ equipped with the distribution $\tD$ induced by the double covering $\tQ  \to Q$.
Let $\Aut(\tQ,\tD)$ be the group of diffeomorphisms of  $\tQ$ leaving   $\tD$ invariant. 
Then we have 

\begin{theorem}  The connected component of the identity in $\Aut(\tQ,\tD)$ for radius ratio 3:1 or 1:3 is isomorphic to $G_2$. The $G_2$ action on $\tQ$ 
does not descend to $Q$, but its   restriction to the maximal compact $K \subset G_2$ does, covering the   $\SO_3 \times \SO_3$ action on $Q$. For any other radius ratio (other then 1:1) $\Aut(\tQ,\tD)$ is ismorphic to $K$.  
\end{theorem}

This theorem was communicated to us by  Robert Bryant, for whom it is but a  variation on a theme of E. Cartan's work on the method of equivalence, contained in his notoriously difficult  ``Five Variables Paper" \cite{C1}  from 1910 .   Bryant  wrote to us recently: 
\begin{quote}
``Cartan himself  gave a geometric description of the flat
$G_2$-structure as the differential system that describes space
curves of constant torsion $2$ or $1/2$ in the standard unit 3-sphere.
(See the concluding remarks of Section 53 in Paragraph XI in the
Five Variables Paper.)  One can easily transform the rolling balls
problem (for arbitrary ratios of radii) into the problem of curves
in the 3-sphere of constant torsion and, in this guise, one can
recover the 3:1 or 1:3 ratio as Cartan's torsion 2 or 1/2 with
a minimum of fuss.  Thus, one could say that Cartan's calculation
essentially covers the rolling ball case.'' 
\end{quote}

\bn\centerline{$*\qquad*\qquad*$}

\bn 

Our main purpose  in this note is to try to explain  this beautiful and mysterious theorem 
in a direct manner which does not appeal  to   Cartan's method of
equivalence. We consider it an  expansion  of Section 4 in  
Bryant's lecture notes \cite{Bryant}.  Our contribution consists basically of a  description of two constructions of  $(\tQ, \tD)$ with 
a built-in  $G_2$-invariance.   Using these constructions we show here   that, for radius ratio 3:1 or 1:3,  $G_2$ is {\em contained} in $Aut(\tilde Q, \tilde D)$.  But
we do not know how to show, without the more sophisticated  Cartan's methods 
(or its variants such as those of Tanaka)
 that  $G_2$ is the   {\em full} identity component
of $Aut(\tilde Q, \tilde D)$, 
nor that for radius ratio different from 3:1, 1:3 or 1:1, $Aut(\tilde Q, \tilde D)$ is not larger than $K$.

A secondary purpose  of this article is to correct an error appearing in the book
\cite{Montbook} by one of us. We had mistakenly said there  that   the symmetry group for the rolling distribution for a ball on a plane   (ratio $1:\infty$) was $G_2$.  

A tertiary purpose is to obtain a bit of a feel for the simplest
 exceptional Lie algebra $\g_2$ and  its Lie groups, and 
 to provide  a refresher course on   roots and weights.

\bn {\bf Structure of Paper.} 
In the next section (section 1) we describe the background and wider context  of the problem,
with references to the literature.   In section 2 we give a detailed   description 
of the distributions associated with the rolling of balls, noting their $\SO_3\times\SO_3$-symmetries. 
In section 3 we describe  the homogeneous distributions of a Lie group $G$
in terms of data   $(G, H, W)$, where $H\subset G$ is a closed subgroup and $W\subset\g/\h$ is an $H$-invariant subspace. 
We then  identify this data for the rolling distribution with respect  the group $G=\SO_3\times\SO_3$.  In section 4 we use the root diagram of $G_2$ to give our first construction of a $G_2$-invariant distribution data $(G_2, P, W)$. The identification of the resulting $G_2$-homogenous distribution on $G_2/P$ with $(\tilde Q, \tilde D)$
 amounts to the embedding of  
 $\so_3\times\so_3$ in $\g_2$ and is   the subject of section 5 (and Appendix B) which forms the heart
 of this article.   
In section 6 we give the second $G_2$-invariant construction of $(\tilde Q, \tilde D)$ an
explicit construction applying   projective geometry to the space of
purely imaginary  split octonions $V$, the lowest dimensional non-trivial representation space
for $G_2$.  Appendix C is historical.  Following suggestions by Bryant
we looked into Cartan's thesis and found that much of content of section 6,
and hence of the rolling distribution already appears there.

\bn\centerline{$*\qquad*\qquad*$}

\bn

Despite all our efforts,  the   ``$3$'' of the ratio $1:3$ remains    
  mysterious.  In this article it  simply arises  out of the structure constants for 
   $G_2$ and   appears  in the construction
of the embedding of $\so_3\times\so_3$ into $\g_2$ (section 5 and Appendix B).
 Algebraically speaking, this  `3'   traces   back   to the 3 edges
 in $\g_2$'s  Dynkin diagram
 and the consequent relative positions of   the long and short
 roots in the root diagram (see figure 2 below) for  $\g_2$ which the Dynkin diagram is encoding.   

 \bn {\bf Open problem.} Find a geometric or dynamical 
interpretation for the ``$3$'' of the $3:1$ ratio.

\bn 

For work in this direction see Agrachev \cite{Agrachev} and also Kaplan and Levstein
\cite{Kaplan}.

 \bn {\bf Acknowledgements.} 
  Robert 
Bryant has   been crucial, at various key steps along the way,
in steering us in the right direction.  Marti Weissmann
supplied us with    key information regarding $G_2$, and the 
crucial Vogan reference.  

 \section{History and Background}

{\bf On distributions.}  By a distribution we mean here a linear subbundle of 
the tangent bundle of a manifold. 
The distributions    first encountered    are  usually
the integrable and the contact distributions and  have infinite dimensional symmetry groups.   In dimension 
$5$  we first encounter distributions whose symmetry groups  
 are  finite-dimensional. Indeed, the generic  distribution of rank $2$ or $3$  in $5$ dimension  has  {\em no} local symmetries. 
 Cartan  \cite{C1} investigated  rank $2$ and
 $3$ distributions in $5$ dimensions in   detail.  
 The growth vector of 
  a generic rank $2$ distribution on a $5$-dimensional manifold,  at a generic
  point of that manifold,   is $(2, 3, 5)$.  
  This  ``$(2,3,5)$  at a point''
    means that   if $X, Y$ are any local vector fields spanning
  the distribution in a neighbhorhood of the point , then $[X,Y] = Z$
  is pointwise linearly independent of $X, Y$ (in a neighborhood of  the point) 
  and $X, Y, Z, [X, Z], [Y, Z]$ span the tangent bundle in a neighborhood of  the point.
  Cartan  worked out the  
complete   local  invariants -- analogues of the Riemann curvature
tensor -- for these $(2,3,5)$ distributions. 
For the distribution's symmetry group to act transitively   all of
Cartan's  invariants must be constant.  To get the maximal dimensional  symmetry group 
all Cartan's  invariants must vanish, in which case we call 
the distribution ``flat''.   Any such distribution
is locally   diffeomorphic to  that of the  ``Carnot group'' distribution
associated to the unique graded nilpotent Lie group $\mathfrak{n} = \mathfrak{n}_{2,3,5}$ of this same  growth, and  its  local  symmetry
algebra   is $\g_2$.
(By  the ``local symmetry algebra'' of a distribution we mean the algebra of vector fields
$X$ satisfying $[X, \Gamma(D)] \subset \Gamma(D)$ where $\Gamma(D)$ is the sheaf of local sections of vector fields tangent to the distribution.)
 
As mentioned in the above quote from Bryant,  Cartan \cite{C1} 
presented several geometric realizations of the flat case.   Bryant  and Hsu 
\cite{BryantHsu}  (see section 3.4) pointed out the
 rolling incarnation of $G_2$.  A $(2,3,5)$
 distribution will    arise whenever  one  rolls one Riemannian surface on  another 
 provided their Gaussian curvatures are not equal.    
 The Cartan invariants vanish if and only if  the ratio of their    curvatures are  $1:9$ .
  Hence the $1:3$ radii for spheres. We could also achieve the maximal
  local symmetry algebra $\g_2$ by 
 rolling two hyperbolic planes along  each other, provided their ``radii''
 are in the ratio $i: 3i$.    
More history, and more instances of the flat $G_2$ system
are  explained in Byrant \cite{Bryant}.

Non-integrable rank $2$ distributions in dimension $n$ ($n > 3$)  admit special families
of integral curves known as ``singular'' or ``abnormal'')
(\cite{Montbook}).  These are curves which
admit no local variations through integral curves and having endpoints fixed.    In the case of $(2,3,5)$ distributions
there is precisely $1$ singular curve (up to reparameterization)
through every point  {\it  in every direction} tangent to $D$. In the case of
rolling one Riemannian surface along another, these singular curves correspond
to rolling along geodesics.   Using  the symplectic geometry associated to
variations of    singular curves
 Zelenko and Agrachev have been  able to 
 rederive  Cartan's (2,3,5) invariants. See   \cite{Agrachev} and references therein. 

Tanaka and his school have established a wonderful
generalization of the passage from the flat nilpotent model $\mathfrak{n}_{2,3,5}$
to $\g_2$.      
Associated to each point $p$ of a manifold endowed with a non-integrable
distribution there is a graded  nilpotent Lie algebra $\m = \m (p)$ called  the
`nilpotentization' of the distribution, or sometimes the ``symbol algebra''.  The dimension of $\m$ is that of the underlying manifold.   
Call the distribution ``of type $\m$'' if the different algebras $\m(p)$
are all isomorphic to the same $\m$, i.e. the isomorphism type does not change
from point to point.  
(Every   (2,3,5) distribution is of type $\mathfrak{n}_{2,3,5}$.)   Associated to each
graded nilpotent $\m$ there is graded Lie algebra $\g \supset \m$, possibly infinite dimensional, 
called the `prolongation' of $\m$ and built from  $\m$ in  a purely algebraic manner.
This  $\g$ represents, roughly speaking,  
the maximal possible symmetry  of a distribution of type $\m$:
every symmetry algebra for a type $\m$-distribution, after applying a 
grading process to it,  
must be a subalgebra of $\g$.  The prolongation of
the (2,3,5) algebra is $\g_2$, and this fact can be viewed 
as the algebraic restatement of Cartan's work on the flat model. 
This Tanaka prolongation method thus yields a proof that $Aut(\tilde Q, \tilde D)
\subset G_2$ in theorem 1, alternative to Cartan's proof. 
  Yamaguchi \cite{Yamaguchi} has classified all $\m$'s whose $\g$'s are simple. 
  To each of these pairs $(\m, \g)$ is associated an intricate differential geometry
  and most of these have not been explored in any detail.

\bn{\bf On $G_2$.}  The Lie algebra $\g_2$  is
the smallest   of the  exceptional simple Lie algebras.
In 1894 Killing uncovered
the   existence of the root lattice for   $\g_2$'s,
but without establishing the existence of the
corresponding Lie algebra.   Cartan, in his thesis,   established the  existence  of $\g_2$
in one page of  his thesis  \cite{C2}.  He did so by constructing
the $7$-dimensional representation of $\g_2$, in a way which is   closely related to our second ``projective split octonion'' 
model for $\tilde Q$, the universal cover of the rolling space.
We have devoted   appendix C to this page of his thesis
and its connection with this  second model.  
  In 1914
Cartan \cite{C3} showed that $G_2$ can be realized as the automorphism
group of the octonions.  For our split $G_2$ he used `split
octonions'.    The compact form of $G_2$ 
  appears  in the Berger list of
potential holonomy groups of Riemannian metrics. Recently,  the compact  $G_2$ has been 
featured in   string theories, but perhaps that fad has passed already.

\section{Distribution for rolling balls}

\subsection{The distribution}

 \bn

Take the first  ball to be stationary, of radius $R$, with its center at the
origin. Roll  a second ball  of radius $r$ on the first ball. The position of the
second ball is given by an isometry (rigid motion)
$\varphi_{(g,\x)}:\R^3\to\R^3$, mapping a point $\P$ to
$$\bp = \varphi_{(g,\x)}(\P)=g\P+(R+r)\x,$$ where $(g,\x)\in\SO_3\times
S^2$.  Here, $R\x$ is the point of contact of the two balls,
$(R+r)\x$ is the center of the second ball and $g\in\SO_3$ describes the
rotation of the second ball relative to its initial position.  See figure 1 in the introduction. Thus the configuration
space $Q$ for our rolling problem has been identified with the manifold
$\SO_3\times S^2$.  (For a visceral account of rolling a sphere on a plane,
accessible to upper division undergraduates,  we recommend
\cite{Hammersley}.)


\bn

Let $(g_t,\x_t)\in Q$ be a  differentiable rolling motion.  Let $\o_t\in
\R^3\cong\so_3$ be the angular velocity of the rolling ball relative to
its center, measured with respect to 
inertial axes.  In other words,  if $\P$  is a material point fixed 
on the second ball, $\dot\P=0$, and if
we write  $\bp_t=g_t\P$,  then $\dot \bp=\dot gg^{-1}\bp=\omega\times
\bp$.  Then we have

\begin{proposition} Let $Q=\SO_3\times S^2$ be the configuration space of two rolling balls of 
radii $R$ and $r$. Let $\rho=R/r$ . 
Then a curve  $(g_t,\x_t)\in Q$ describes a rolling motion without slipping and spinning iff
\mn \begin{quote}
(1) $(\rho+1)\dot \x=\omega\times \x$ (no-slip condition),

(2) $\langle\omega,\x\rangle=0$ (no-spin condition,
 i.e. $\omega$ need to be tangent  to the stationary ball at $R\x$).
\end{quote}
\end{proposition}

\pf (1) The contact point between the two balls is $\bp=R\x$ on the
first ball, $\P=-g^{-1}r\x$ with respect to the second ball. For
non-slip, their velocities must match: $\dot\bp= g\dot\P.$ Now
$\dot\bp=R\dot\x$ and
$$\dot\P= [-{d\over dt}g^{-1}]r\x-g^{-1}r\dot\x=g^{-1}\dot gg^{-1}r\x-g^{-1}r\dot\x=g^{-1}r(\o\X\x-\dot\x),$$
hence the non-slip condition $\dot\bp= g\dot\P$ is equivalent to $R\dot\x=r(\o\X\x-\dot\x),
$
from which (1) follows.

\mn (2) Let $\bf P$ be a material point {\em fixed} on the second ball
($\dot\P=0$). From the inertial point of view, which is to say,
from the point of view  of the first ball with origin at its center,
the position of this material point is  $\bp=g\P+(R+r)\x$, and
so its  velocity
$$\dot\bp=\dot g\P+(R+r)\dot\x=\dot g
g^{-1}[\bp-(R+r)\x]+(R+r)\dot\x=\o\X[\bp-(R+r)\x]+(R+r)\dot\x.$$ Using
the no-slip equation, $(R+r)\dot \x=r\omega\times \x$, we get
$$\dot\bp=\o\X[\bp-(R+r)\x]+r\o\X\x=\o\X(\bp-R\x).$$ 
The equation
$\dot\bp=\o\X(\bp-R\x)$ asserts  that the instantaneous motion of the second ball is
 a rotation whose
axis of rotation (a line) passes through the point of contact $R\x$,
in the direction of $\o$ and with angular velocity of magnitude
$\|\o\|$. The no-spin condition is 
that the second ball does not spin about
the   point of contact of the two balls, which is to say that $\o$ should have no
component orthogonal to the common tangent plane of the two balls,
i.e. $\langle\omega,\x\rangle=0$. \qed

The two conditions in the last Proposition define together a rank 2 distribution on $Q$. This is the rolling distribution.

\vskip .3cm

\subsection{The ``obvious'' symmetry}

The group $\SO_3\times\SO_3$ acts  on $Q$ by
$\varphi_{(g,\x)}\mapsto g'\circ\varphi_{(g,\x)}\circ g''^{-1},$ where
$g',g''\in\SO_3.$ In terms of $(g,\x)$ this 
action is
$$(g,\x)\mapsto (g' g g''^{-1},g'\x), \quad g',g''\in\SO_3.$$
This action is transitive and preserves the rolling distribution $D$ for any value of $\rho=R/r$.  The proofs of these assertions are easy and left as exercises.

\section{Group theoretic description of the rolling distribution
}

\bn

In the previous section we wrote down a distribution $D$ on $Q=\SO_3\times S^2$, depending on the real parameter $\rho$.  We showed
that $Q$ admits   an $\SO_3\times\SO_3$-transitive action which preserves the distribution. Our aim in this paper is to show that for two specific values of the parameter, $\rho=3$ and $\rho=1/3$, the distribution admits a larger {\em local} group of symmetries, namely the group $G_2$.  We do so by defining a $G_2$-homogeneous space $\tQ=G_2/P$, together with a $G_2$-invariant rank 2 distribution $\tD$ on it. 
We then define a 2:1 covering map $\tQ\to Q$ which maps $\tD$ to $D$. Furthermore, the group $G_2$ contains a maximal compact  subgroup $K\subset G_2$ which is a double cover of $\SO_3\times\SO_3$, such that the map $\tQ\to Q$ is $K$-equivariant  (with respect to the covering homomorphism $K\to\SO_3\times\SO_3$). 
The  constructions 
are most easily done on the group level or on the Lie algebra level. We describe in what follows the general set up required for ``working on the group level'' and 
then calculate the group theoretic data corresponding to the rolling distribution. 

\bn 

Let $G$ be a Lie group. A ``$G$-homogeneous distribution'' is a pair  $(Q,D)$ where  $Q$  is a manifold on which $G$ acts  transitively  and  $D\subset TQ$ is  a $G$-invariant distribution. Fixing  a base point 
 $q_0\in Q$ with isotropy  $H\subset G$ we obtain a $G-$equivariant identification    $ Q\cong G/H$  and an $H$-equivariant identification $T_{q_0}Q \cong \g/\h$
 where $\h \subset \g$ denote   the Lie algebras corresponding to $H \subset G$.
 Then  $D_{q_0}\subset T_{q_0}Q$ corresponds to an $H$-invariant subspace $W\subset\g/\h$. 
 In this way  every  $G$-homogeneous distribution $(Q,D)$ corresponds to data 
 $(G,H,W)$, where  
$H\subset G$ is a closed subgroup with Lie algebra $\h\subset\g$ and 
$W\subset\g/\h$ is an $H$-invariant subspace.
 The adjoint action of $G$ defines an equivalence relation on the set of pairs  $(H,W)$ so that 
different choices of base points on $Q$ correspond to equivalent pairs  $(H,W)\sim(H',W')$. 
 Conversely,  given the data $(G,H,W)$, we can construct a $G$-homogeneous distribution $(Q,D)$ by letting $G$  act by left translations on the
 right $H$-coset space $Q:=G/H$, and define a $G$-invariant distribution
 $D\subset TQ$  using the $G$-action to  push   $D_{[e]}:=W\subset\g/\h\cong
 T_{[e]}(G/H)$ around to all other points of $Q$.

On the level of  Lie algebras, the  data $(\g,\h,W)$ determines $(Q,D)$
  up to a cover.  If, as in our case of $\g= \so_3 \oplus \so_3$, the simply
  connected Lie group $G$ realizing $\g$ is compact,
  then there are only finitely many homogeneous distributions    $(G, H, W)$ which realize the   given  Lie algebraic data $(\g, \h, W)$.

\bn 
  
We now  determine  the data $(G,H,W)$ corresponding to the  rolling  distribution
$(Q,D)$ 
of section 2.1.   Here $G=\SO_3\X\SO_3$, $Q=\SO_3\times S^2$, $\dim H=1$, $\dim W =2$. 
Identify   the Lie algebra of $\SO_3\times\SO_3$ with $\R^3\times\R^3$, the set of pairs of angular velocities $(\omega', \omega'')$, with 
Lie bracket given by the cross product:
$$[(\omega', \omega''),(\eta', \eta'')]=(\omega'\X\eta',
\omega''\X\eta'').$$  The first factor $\omega '$ corresponds to the
first (stationary) sphere, of radius $R$, while the second $\omega ''$ factor corresponds to second (rolling) sphere of   radius $r$. 

Fix a base point, say $(1,\e_3)\in\SO_3\times S^2=Q$.  The
isotropy at this base point is the circle  subgroup $H$ consisting of elements of the form $(h,h)$, where $h$ is a
rotation around the $\e_3$ axis, so $\h=\R(\e_3,
\e_3)\subset\R^3\X\R^3.$ Using the Killing metric on $\g=\so_3\X\so_3=\R^3\X\R^3$ we can
identify $\g/\h\cong\h^\perp$, so that the plane of the distribution at the
base point is given by some 2-plane in $\h^\perp.$ Let us
determine explicitly this 2-plane.

\begin{proposition} The rolling distribution  on $\SO_3\times S^2$ 
corresponding to rolling a ball of radius $r$ along one of radius $R$  is given by the 2-plane
in $\R^3\times\R^3$ (the Lie
algebra of $\SO_3\X\SO_3$) defined by the equations
$$\langle\omega',\e_3\rangle=\langle\omega'',\e_3\rangle=0, \quad
\rho \omega'+\omega''=0,$$
where $\rho = R/r$.
\end{proposition}

\pf Since $\h\subset\R^3\times\R^3$ is  generated by the vector $(\omega', \omega'')=(\e_3,\e_3)$ and the Killing metric corresponds to some multiple of the standard metric on $\R^3\times\R^3$, $\h^\perp\subset \R^3\times\R^3$  is given by the equation $\langle\omega'+\omega'',\e_3\rangle=0.$

From the formula for the $\SO_3\X\SO_3$-action in \S2.2 we get the
infinitesimal action at the base point $$\omega=\omega'-\omega'',
\quad\dot \x=\omega'\X \e_3.$$ Substituting these into the
rolling  conditions at the base point (see \S2.1),
$$\langle\omega,\e_3\rangle=0, \quad (R+r)\dot \x= r\omega\times \e_3,$$
we obtain
$$\langle\omega'-\omega'',\e_3\rangle=0, \quad  [R\omega'+ r\omega'']\times \e_3=0.$$
Adding the condition of orthogonality
to $\h$, $\langle\omega'+\omega'',\e_3\rangle=0,$
we obtain the above equations.$\Box$

\bn

\subsection{Shrinking the group.} 
\label{shrinking} The following observation will be key to proving
that part of   theorem 1 which we are going to prove, namely
that $G_2 \subset Aut(\tilde Q, \tilde D)$ for $\rho=3$ or $1/3$.
Suppose that   $(Q,D)$ is a $G$-homogeneous distribution with $G$-data $(H,W)$. 
Let $G_1\subset G$ be a subgroup for which the restriction of the $G$-action on $Q$ to $G_1$ is still transitive. 
Then $(Q,D)$ is also $G_1$-homogeneous distribution and its   $G_1$-data is $(H_1,W_1)$  where $H_1=H\cap G_1$ and $W_1\subset \g_1/\h_1$
corresponds to 
$W$ under the  linear isomorphism   $\g_1/\h_1\to\g/\h$
induced by the diffeomorphism $Q = G/H = G_1/ H_1$.
Since $(Q, D)$ has not been changed, it follows that
the $G$-data $(H, W)$ and the $G_1$-data $(H_1, W_1)$
yield diffeomorphic manifolds with distributions. 
 At the Lie algebra level, this discussion asserts that
$(\g_1, \h_1, W_1)$ and $(\g, \h, W)$ define manifolds-with-distributions
which are diffeomorphic up to a cover.  To prove that
$G_2 \subset Aut(\tilde Q, \tilde D)$ we  will be applying  this observation to 
the case $\g_1 = \so_3 \oplus \so_3 \subset \g_2$.

\section{A $G_2$-homogeneous distribution }

\bn

We now describe the other main actor  in this paper,
a distribution with Lie algebraic data $(\g_2, \p, W)$. 
Please see  the root diagram of $\g_2$ in figure 2.  This diagram will be explained
immediately below.  The decorations on the diagram are used to  indicate
the Lie algebraic data and  will   be explained a bit later. 

\bn

\begin{center}
\unitlength=1cm
\begin{picture}(5,4)(-2,-2)
\put(0,0){\circle*{.4}}

\put(.4,.75){$\oplus$}
\put(-.5,-.866){\circle*{.25}}
\put(.6,1.2){$\sigma_1$}
\put(-1.2,-1.3){$-\sigma_1$}

\put(-.5,.866){\circle*{.25}}
\put(.4,-.95){$\oplus$}
\put(-.95,1.2){$\sigma_2$}
\put(.3,-1.3){$-\sigma_2$}

\put(-1,0){\circle*{.25}}
\put(1,0){\circle{.25}}
\put(1.3,-0.1){$\sigma_3$}
\put(-2,-0.1){$-\sigma_3$}

\put(1.4,-.866){\circle{.25}}
\put(-1.5,.866){\circle*{.25}}
\put(1.6,  -1.1){$-\lambda_1$}
\put(-2.2,  0.866){$\lambda_1$}

\put(1.4,.866){\circle{.25}}
\put(-1.5,-.866){\circle*{.25}}
\put(1.7,.9){$\lambda_2$}
\put(-2.4,-1.1){$-\lambda_2$}

\put(0,1.5){\circle*{.25}}
\put(0,-1.5){\circle*{.25}}
\put(-.2,1.8){$\lambda_3$}
\put(-0.5,-2.1){$-\lambda_3$}

\put(0,-1.4){\line(0,1){3}}
\put(-1.5,  -0.866){\line(5,3){2.77}}
\put(-1.5,  0.866){\line(5,-3){2.77}}

\put(  -1, 0){\line(1  , 0   ){1.85}}
\put(-.5,-.866){\line(3,5){.96}}
\put(-.5,.866){\line(3,-5){.96}}
\end{picture}

\bn

{\small Figure 2: The root diagram of $\g_2$.}
\end{center}

\bn{\bf A reminder of the meaning of the root diagram.} The plane in which the  diagram is drawn is  the dual of a Cartan subalgebra $\t\subset\g_2$. A Cartan
subalgebra of a semi-simple Lie algebra $\g$ is a maximal abelian
subalgebra $\t\subset\g$ of semi-simple elements, i.e. each
$ad(T)\in\End(\g)$, $T\in\t$, is diagonalizable.  In the case of
$\g=\g_2$, $\t$ is 2-dimensional, hence the subscript 2 in $G_2$, the
{\em rank} of the group.  The root diagram of $\g$ encodes the adjoint
action of $\t$ on $\g$, from which one can recover the whole structure
of $\g$.

The commutativity of the Cartan
subalgebra $\t$ implies that the diagonalizable endomorphisms
$ad(T)\in\End(\g)$, $T\in\t$, are {\em simultaneously} diagonalizable,
resulting in a $\t$-invariant decomposition
$$\g=\t\oplus\sum_\alpha\g_\alpha,$$
where each  $\g_\alpha\subset\g$ is a 1-dimensional subspace  of $\t$-common eigenvectors called a {\em root space}.
The corresponding eigenvalue  depends linearly on
the acting element of $\t$, so is given by a linear functional $\alpha\in\t^*$, called {\em root}.  Thus
$$[T,X]=\alpha(T)X, \quad T\in\t,\quad X\in\g_\alpha.$$

 When we draw the root diagram in $\t^*$
we use the Killing metric in $\g$ to determine the size of the roots
and especially  the angles between them. The Killing metric in $\g$ is the inner
product $\langle X,Y\rangle=-{\rm tr}(ad(X)ad(Y))$.  It is
non-degenerate (this is equivalent to semi-simplicity) and its
restriction to $\t$ is positive definite.

\bn{\bf Example of $\g=\sl_3(\R).$} The more familiar example of
$\sl_3 (\R)$ is  useful   to keep
in mind before   proceeding with $\g_2$.   The Lie algebra $\sl_3(\R)$ is the 
vector space of 
of 3 by 3 traceless real matrices with Lie bracket the usual matrix Lie bracket.
It is  Lie algebra of the Lie group $\SL_3(\R)$ of 3 by 3 real matrices with determinant $1$.   Like $\g_2$, the Lie
algebra $\sl_3(\R)$ is a non-compact split form of its
complexification ($\sl_3(\C)$) and has  rank 2. We take as a Cartan
subalgebra the subspace  $\t\subset\sl_3(\R)$ of traceless diagonal
matrices,
$$\t:=\{\left(\begin{array}{ccc} t_1&0&0\\0& t_2&0\\0&0& t_3\end{array}\right) |   t_1+ t_2+ t_3=0, t_i\in\R\}.$$

$\sl_3(\R)$ has    6 roots: 
$$\alpha_{ij}:= t_i- t_j\in\t^*, \quad i \neq j, \quad
i,j\in\{1,2,3\},$$ with corresponding root spaces
$$\g_{\alpha_{ij}}=\R E_{ij}, $$ where $E_{ij}$ is the matrix whose
$ij$ entry is 1 and all of whose other   entries are 0. The corresponding root space
decomposition $$\sl_3=\t\oplus\sum_{i\neq j}\g_{\alpha_{ij}},$$ is
just the decomposition of a matrix as a diagonal matrix plus its off
diagonal terms.  The metric induced on $\t$ by the Killing metric is
some multiple of the standard euclidean metric, so that $\langle T,
T'\rangle=c\sum_i t_i t'_i$ for some $c>0$.

\bn

\begin{center}
\unitlength=1cm
\begin{picture}(5,4)(-2,-2)
\put(0,0){\circle*{.4}}

\put(-1.5,.866){\circle*{.25}}
\put(1.3,-1){$\oplus$}
\put(-2.3,1){$\alpha_{12}$}
\put(1.2,-1.3){$-\alpha_{12}$}

\put(0,1.5){\circle*{.25}}
\put(0,-1.5){\circle{.25}}
\put(-.3,1.8){$\alpha_{13}$}
\put(-.5,-2){$-\alpha_{13}$}

\put(1.4,.866){\circle*{.25}}
\put(-1.7,-1){$\oplus$}
\put(1.6,1){$\alpha_{23}$}
\put(-2.2,-1.3){$-\alpha_{23}$}

\put(-1.4,  -0.8){\line(5,3){2.9}}
\put(0,-1.4){\line(0,1){3}}
\put(-1.5,  0.866){\line(5,-3){2.8}}
\end{picture}

\bn

{\small Figure 3: The root diagram of $\sl_3$}

\end{center}

\bn

\bn{\bf Reading the root diagram}.  One can read much of the structure
of $\g$ from its root diagram  in a formula-free manner. Here
is the key observation. Let $\alpha, \beta$ be two roots with
(non-zero) root vectors $E_\alpha\in\g_\alpha,$ $ E_\beta\in\g_\beta$.
That is, $$[T,E_\alpha]=\alpha(T)E_\alpha,\quad T\in\t, $$ and
similarly for $\beta$. It then follows immediately from the Jacobi
identity that
 $$[T,[E_\alpha, E_\beta]]=(\alpha+\beta)(T)[E_\alpha, E_\beta].$$ This means that
 \begin{itemize}

 \item[(1)] if  $\alpha+\beta\neq 0$ and is not a root then $[E_\alpha, E_\beta]=0$;

 \item[(2)] if  $\alpha+\beta\neq 0$ and is a root then $[E_\alpha, E_\beta]\in\g_{\alpha+\beta}$;

 \item[(3)] if $\alpha+\beta= 0$, i.e. $\beta=-\alpha$, then $[E_\alpha, E_\beta]\in\t$.
 \end{itemize}

This set of 3 conclusions permit us to see at a glance from the
diagram a fair amount of the structure of $\g$.  In the last two cases
one can further show that $[E_\alpha, E_\beta]$ is non-zero and
determine, with some calculations, the actual bracket, as will be
illustrated in Appendix B.

\bn{\bf Example: reading the root diagram of $\sl_3$.} Let us consider
the subspace $\p\subset\sl_3$ spanned by $\t$ and the root spaces
corresponding to the roots marked with dark dots in figure 3.

 The diagram shows that $\p$ is a 5-dimensional subalgebra, i.e. it is
closed under the Lie bracket (there are 4 dark dots, but remember that
the thick dot at the origin stands for the 2-dimensional Cartan
subalgebra). Indeed, $\p$ is the subalgebra of upper triangular
matrices (including diagonal ones), with corresponding subgroup
$P\subset\SL_3$, the subgroup of upper triangular matrices with
determinant=1.  The quotient space $\SL_3(\R)/P$   can be
identified with the space $F$ of full  flags in $\R^3$.  A full
flag is a  pairs
$(l,\pi)$, where $l$ is a line and $\pi$ is a plane, and
$l\subset\pi\subset\R^3$.    The ``standard flag'' consisting
of the $x$ axis sitting inside the $xy$ plane has isotropy group $P$. The tangent space to $F$ at this  base
point  is naturally identified with $\sl_3/\p$,
represented in the root diagram by the remaining three light dots.
Two of the light dots are marked +. The diagram shows that the root spaces corresponding to these roots
span a $\p$-invariant 2-dimensional subspace of $\sl_3/\p$ which Lie
generates the root space associated with the third light
dot. 
This means that we have on $F$ an
$SL_3(\R)$-invariant rank 2 contact distribution, i.e. a
non-integrable distribution that Lie generates the tangent bundle. 

This distribution can be geometrically interpreted
as  the  ``tautological"
contact distribution on $F$.  This distribution  is spanned by two vector fields,
corresponding to the two +s in figure 3.  One vector field generates the flow
in which the line   $l$ spins within   the plane $\pi$, while  the plane remains  fixed.  The other vector field generates the flow in which   the plane $\pi$ rotates about
the line $l$, while the line remains fixed.

\bn{\bf Reading the $\g_2$ diagram.}  
Now let us draw conclusions in a
similar fashion from the $\g_2$ diagram.  There are   twelve roots 
in the diagram (figure 2) and so 12 root spaces.  The rank of $\g_2$ is $2$
and so the dimension of $\g_2$ is $14 = 2 + 12$.  
  Consider the 9-dimensional
subspace $\p\subset\g_2$ spanned by $\t$ and the root spaces
associated with the roots marked by the dark dots in the diagram of figure 2.
Then the diagram shows that
\begin{itemize}
\item $\p$ is closed under the Lie bracket, i.e. is a subalgebra (a
so-called parabolic subalgebra, a subalgebra containing a Cartan
subalgebra).
\item Let $P\subset G_2$ be
the corresponding subgroup. It follows that $G_2$ has a 5-dim
homogeneous space $G_2/P$, whose tangent space $\g_2/\p$ at a point is
represented by the remaining 5 light dots.

\item Two of the light
dots are marked with +. The diagram shows that their root spaces generate a 2-dim $\p$-invariant
subspace  $W_1\subset\g_2/\p$, hence a $G_2$-invariant rank 2 distribution on $G_2/P$.

\item This
distribution is not integrable, in fact, it is a distribution of type $(2,3,5)$,
since the diagram shows that bracketing once gives the light dot marked with
$\sigma_3$ and bracketing again gives the remaining two light dots.

\end{itemize}

To summarize, we have  assembled the ingredients for the data
$(G_2, P, W_1)$.   To prove the theorem
is to provide the geometric interpretation
of this distribution as $(\tilde Q, \tilde D)$ of theorem 1. 
The first step in doing so is to embed $\so_3\oplus\so_3$ in $\g_2$.

\section{The  maximal compact subgroup of $G_2$}

\subsection{Algebraic strategy of the proof.}

In the previous sections we   assembled the  Lie algebraic   data,
$(\so_3 \oplus \so_3, \h, D)$ and $(\g_2, \p, W)$
with corresponding group data $(K, H, D)$ and   $(G_2, P, W)$.  
The key to  theorem 1  is to embed
$\so_3\oplus\so_3$   in $ \g_2$.  This embedding is constructed
in   the next section.    Once established, we obtain
a diffeomorphism between corresponding distributions
by following the   observation made in  section \ref{shrinking}.

We recap that observation,  adding a bit of topology.  
Suppose that $\k \subset \g$
and $\p \subset \g$ are   Lie subalgebras of the Lie algebra $\g$.
Suppose that the natural map  $\k/\k \cap \p \to  \g/ \p$ is a linear 
isomorphism.  Let  $W \subset \g/\p$ be an $ad$- $\p$-invariant subspace 
and   $W_1 \subset \k/\k \cap \p$
 the corresponding subspace.
Then   we will say that the   Lie algebraic distributional data
$(\g, \p, W)$ and
$(\k, \k \cap \p,  W_1)$ are isomorphic.  If the corresponding 
connected Lie groups are $K \subset G$ and
{\it if } $K$ is compact, then we can conclude
that the data $(G; P, W)$ and $(K; K \cap P, W_1)$
define isomorphic manifolds with distributions.
For when $K$ is compact we have that $K/K \cap P$ is a
compact and open submanifold of $G/P$ and hence is diffeomorphic
to $G/P$.  And under this diffeomorphism
the distribution corresponding to $W$ is the same as the one represented
by $W_1$.  

The compactness assumption on $K$ is neccessary to conclude
that $G/P = K/(K \cap P)$. 
Think of the case $K = \C^* \subset G = SL(2,\C)$
where $G$ acts on the sphere $Q = \C \cup \{ \infty \}$  by Mobius transformations
and where $\C^*$ corresponds to the complex scalings $z \mapsto \lambda z$,
$\lambda \ne 0$. 
The fixed points of the $\C^*$-action are  $0, \infty$.
The  $\C^*$ orbit through any   point $z_0 \ne 0, \infty$
is open, being the whole sphere minus the two fixed points.
 Thus $G/P \ne K/ (K \cap P)$ where $P$ is the isotropy group of $z_0$. 
 But we still have $\k/(\k \cap \p) = \g/ \p$ since the orbit of $z_0$ is open. 
 
 Lie algebraic data   defines the corresonding Lie group data  only  up to a covering.
 We can insure that there are only finitely many such coverings
 by knowing that  $\k$, like $\so_3 \oplus \so_3$, is  the compact
real form of its corresponding complex Lie algebra.
For in this case there are only a finite number of 
connected  Lie groups $K$  with    Lie algebra $\k$,
all of these being compact and covered by the simply connected $K$.
To say this in another way, suppose we are given  Lie algebra data
$(\g, \p, W)$ and $(\k, \k \cap \p,  W_1)$ as above, and
suppose that $\k$ is a compact real form.  Let 
$(G, P, W)$ and $(K, H, W_1)$ denote {\it any}   Lie-group data realizing these
 respective  Lie algebraic data  {\it where we are no longer assuming that $K \subset G$}.
 Then the   two manifolds-with-distribution which they stand for
are isomorphic up to a finite cover. 
By this we mean, there is a third  manifold-with-distribution $(X, E)$
and covering maps $\pi_G: X \to G/P$,
$\pi_K:  X \to K/H$ such that $\pi_G ^* W = \pi_K ^* W_1 = E$.
Indeed,  we can take   $X$  to be $\tilde G/ \tilde P$
where $\tilde G$ is the unique simply connected Lie group with algebra
$\g$. 

To establish theorem 1  we will apply these considerations
to the case  $G = G_2$ and   $K \subset G_2$
  its maximal compact subgroup.  
We will show that $(\k, \k \cap \p)$
is isomorphic to $(\so_3\X\so_3,  \h = \R(\e_3,\e_3))$.
And we show that under this isomorphism
$W_1 \subset \k/ \k \cap \p$
corresponds to the rolling distribution 
when the   ratios of the rolling spheres are   $1:3$.

\subsection{ Finding   Maximal compacts.} How can we ``see'' a maximal compact subgroup of $G_2$ tangled within   its root diagram? Let us look back again
at the example of $\SL_3(\R)$. Here the
 maximal compact subgroup is $\SO_3$, with Lie algebra $\so_3$, the set of 3 by 3 antisymmetric matrices.
 These are spanned by the vectors $E_{ij}-E_{ji}$, $i>j$. So we see that corresponding to each pair
of ``antipodal" roots  $\pm\alpha_{ij}$ we have one generator of $\k$,  lying   in the sum of the two corresponding root spaces.

More generally, for the ``split'' real form of any semi-simple Lie
algebra (such as our $\g_2$), the situation is similar: we get the Lie
algebra $\k$ of a maximal compact subgroup $K\subset G$ by taking the
sum of 1-dimensional subspaces, one subspace for each pair of
antipodal roots $\pm\alpha$.  In fact, there is a certain
particulary ``nice" choice of root vectors $E_\alpha\in\g_\alpha$
(sometimes called a ``Weyl basis''), so that the sought-for line is
given by $\R(E_\alpha-E_{-\alpha}),$ as in the $\sl_3$ case.

\bn

In the case of $\g_2$ we thus have that \begin{itemize}
\item  $\k$ is the sum of six 1-dimensional subspaces  $\s_{i},\l_{i}$, $i=1,2,3,$ where
$\s_{i}$ lies in the sum of the root spaces corresponding to $\pm\sigma_i$, and
 $\l_{i}$ lies in the sum of the root spaces corresponding to $\pm\lambda_i$.

\item The isotropy of the $K$-action, $H=K\cap P\subset K$,  is given  in the diagram  by
the vertical segment, $\h=\l_{3}$.

\item The distribution plane $W\subset \k/\h$ is generated by  $\s_{1},\s_{2}$ (mod $\h$).

\end{itemize}

We have thus assembled the required ingredients for a   ``distribution data" $(\k,\h,W)$.

\subsection{$\k\simeq\so_3\oplus\so_3$}

Our task here is to define an isomorphism $\k\simeq\so_3\oplus\so_3$
that maps $(\k,\h,W)$ to the data of \S4 with $\rho=3$ or $1/3$. This entails the
decomposition of $\k$ into the direct sum of two ideals, each 
isomorphic to $\so_3$. It would have been quite nice and simple if the
sought-for decomposition of $\k$ had been the decomposition into
``long" ($\l_{i}$) and   ``short" ($\s_{i}$). But this is not the case. For   the
diagram shows that although the $\l_{i}$ generate an
$\so_3$ subalgebra of $\k$, this subalgebra is not an ideal, so is not one
of the summands in the decomposition. And the $\s_{i}$
do not generate even a subalgebra. So we have to work harder, i.e. write
down the precise commutation relations.

\begin{proposition} There is
 a basis $\{S_i,L_i | i=1,2,3\}$ of $\k$, with $S_i\in \s_{i}$ and $L_i\in \l_{i}$,
 such that
$$[L_i,L_j]=\epsilon_{ijk}L_k,\quad [L_i,S_j]=\epsilon_{ijk}S_k,\quad [S_i,S_j]=\epsilon_{ijk}({3\over 4}L_k-S_k),$$
where  $\epsilon_{ijk}$ is the ``totally antisymmetric tensor on 3 indices" (
$\epsilon_{ijk}=1$
if $ijk$ is a cyclic permutation of 123, -1 if anticyclic permutation, and
0 otherwise).
\end{proposition}

The proof of this proposition is relegated to Appendix B. 
It  consists of  simple but tedious
calculations which we could not ``see'' in the diagram. We tried.  We
were reduced to picking up as nice as possible basis for $\g_2$ and
calculating the corresponding structure constants with the help of  Serre  \cite{Serre}. 

\bn Now set $$\e_i':={3L_i+2S_i\over 4}, \quad \e_i'':={L_i-2S_i\over 4}, \quad i=1,2,3.$$
These 6 vectors form a new basis for $\k$ and satisfy the standard $\so_3\oplus\so_3$ commutation relations
\begin{equation}
[\e_i',\e_j']=\epsilon_{ijk}\e_k', \quad
[\e''_i,\e''_j]=\epsilon_{ijk}\e''_k, \quad [\e_i',\e''_j]=0,
\label{std}
\end{equation}
 thus
establishing the desired Lie algebra isomorphism $\k\simeq\so_3\oplus\so_3$.

\begin{corollary} The map $\k\to\so_3 \oplus\so_3$ defined by
$\e_i'\mapsto (\e_i,0),$ $\e_i''\mapsto (0,\e_i),$ $i=1,2,3,$ is a Lie
algebra isomorphism. It maps $\h=\R L_3$ to $\R(\e_3,\e_3)$ and the
2-plane in $\k$ generated by $S_1,S_2$ to the 2-plane in
$\so_3\X\so_3$ defined in the Proposition of \S4 for $\rho=3$.
Interchanging the summands in $\so_3 \oplus\so_3$, i.e. mapping 
$\e_i'\mapsto (0,\e_i),$ $\e_i''\mapsto (\e_i,0),$ correponds to $\rho=1/3$. 
\end{corollary}

The first assertion is eq (\ref{std}). The second assertion
is easily verified using the last Proposition. We have thus
defined a $G_2$-action on some finite cover of the rolling configuration,
one which preserves the pulled-back distribution when the 
radii of the two balls are in the ratio  $3:1$ or $1:3$. QED

\bn{\bf How we came up with the formulae for $\e'_i,
\e''_i$.} The first thing to observe is that since $L_3$ generates the
isotropy $H=P\cap K$ we should have $L_3=\e_3'+\e_3''.$ Since
everything is symmetric in 1,2,3 we conclude that $L_i=\e_i'+\e_i'',$
$i=1,2,3.$ Next since $S_3$ commutes with $L_3$ we should have
$S_3=a\e_3'+b\e_3''$ for some constants $a,b$, and again by symmetry
$S_i=a\e_i'+b\e_i''$, $i=1,2,3.$ Now by using the sought-after commutations
relations for the $\e'_i, \e''_i$ and the known commutations for  $L_i, S_i$ we get that $a,b$
are roots of the equation $x^2+x-3/4=0,$ i.e. $a=1/2, b=-3/2.$ Hence,
$$L_i=\e_i'+\e_i'',\quad S_i=(\e_i'-3\e_i'')/2,\quad i=1,2,3.$$
Inverting these equations we obtain the above equations for $\e_i', \e_i''$.

\section{Split Octonions and the projective quadric realization of $\tilde Q$}

To show that $G_2 \subset Aut(\tilde Q, \tilde D)$
(theorem 1),   it remains to identify $G_2/P$ with
the $\tilde Q = S^3 \times S^2$ of the theorem   
and to show that the covering map $S^3 \times S^2 \to Q$
corresponds to the identification $G_2 /P = K/H$
composed  with the projection $K/H \to (\pm 1, \pm 1) \backslash K / H$.
In order to do these things  will use the fact, discovered by Cartan \cite{C3}
in 1914, that $G_2$
is the   group of automorphisms
of the ``split octonions'' $\tilde \O$.  We will
follow the treatment in the book  \cite{Harvey}, in  the section   ``The
Cayley-Dickson process" (p.104).  There  further   consequences and motivation  
can be found.    

The split octonions $\tilde \O$ are a real  eight-dimensional
algebra with unit and  which is neither associative nor commutative. 
We identify $\tilde \O$ with     $\H^2$, the 2 dimensional quaternionic vector space.
Its   multiplication law is  
\begin{equation}
(a,b)(c,d)=(ac+\bar d b,da+b\bar c), \hskip .6cm a, b, c , d \in \H
\label{product}
\end{equation}
 The unit  $1 \in  \tilde \O$ is   $(1,0) \in \H^2$

 The  automorphism group of  a real algebra is  $A$
  is defined to be 
  the space of nonzero real invertible linear maps  $g: A \to A$
 satisfying  $g(xy) = g(x) g(y)$ for all $x, y \in A$.
  $G_2$ is the automorphism group of $\tilde \O$.  
    
  The unit  $1$ 
  is automatically invariant under any automorphism of $\tilde \O$,
  so that $\R = \R 1 \subset \tilde \O$ is a $G_2$-invariant subspace.
  This subspace has a $G_2$-invariant complement:  
  $$ \tilde \O = \R 1 \oplus V = Re(\tilde \O) \oplus Im(\tilde \O) $$
  In  quaternionic terms: 
  \begin{equation}
  V=Im\tilde\O=Im\H\oplus\H   \subset \H ^2 .
  \label{V}
  \end{equation}
  To see the invariant nature of $V$, we
  use   the split-octonion conjugation  $x \mapsto \bar x$
  defined by  $x=(a,b)\in\tilde\O \mapsto \bar x=(\bar a, -b)$
  for $x \in \tilde\O$.
   Then $x=Re(x)+Im(x)$,
$Re(x)=(x+\bar x)/2 \in \R 1$, and $Im(x)=(x-\bar x)/2$.
Also $x \bar x = - \langle x, x \rangle 1$
where  $\langle x,y\rangle=Re(x\bar y)$ is an inner product of signature $4,4$ on
$\tilde\O$ which is  invariant under the action of $G_2$.     
$V$ is the orthogonal complement of 
$1 \in \tilde \O$  relative to this inner product.
Alternatively, it  can be shown that $x \in V$ if and only if $x^2 = \langle x, x \rangle 1$
 (see \cite{Harvey}, lemma 6.67), proving the
   $G_2$-invariance of $V$.   
$V$   forms  a 7-dimensional inner product 
space of signature $(3,4)$ relative to the restriction of 
$\langle \cdot , \cdot \rangle$. The $G_2$ action on $V$
leaves this inner product invariant, so that  
$G_2$ is realized as   a subgroup of $SO(3,4)$ through its representation
on $V$. 

The maximal compact of $G_2$
is $K = SO(4) = ( SU(2) \X SU(2) )/ \pm(1,1)$. 
See Appendix B and \cite{Vogan}.  Upon restricting from  $G_2$ 
to $K$, the representation  $V$ decomposes into irreducibles according
to   (\ref{V}).  In other words, thinking of $SU(2)$ as unit quaternions, 
  $(q_1, q_2) \in SU(2) \times SU(2) = \tilde K$ (the universal cover of $K$)
   and $(a,b) \in Im(\H) \oplus \H = V$ we
have $(q_1, q_2) \cdot (a,b) = (q_1 a \bar q_1, q_1 b \bar q_2)$.  

In quaternionic terms (\ref{V}) the quadratic form associated to 
our $(3,4)$ inner product on $V$ is 
$$\langle (v, q), (v, q) \rangle = -|v|^2 + |q|^2.$$
Since $K$ acts transitively on the product of spheres
$S^2 \times S^3 \subset Im (\H) \oplus \H = \tilde \O$
we have that 
$G_2$ acts transitively on the null cone
$\{ x = (v, h): \langle x,x\rangle=0, x\neq 0 \}$.
(To see that we can change the `length'
of an $x$ in the null cone using $G_2$,
use the fact that each such  null vector is a nonzero weight vector relative
to some choice of maximal Cartan $T \subset G_2$.
This maximal Cartan then acts on $x$ by scaling.  See the description
following eq. (\ref{basis}) below.) 
Thus $G_2$ acts  
 transitively on the space  of null rays
 $$ C=\{\R^+ x\subset V|\langle x,x\rangle=0, x\neq
0\}\subset P^+ (V) := \hbox{ rays in } V $$. 
This $C$ is a nondegenerate 5-dimensional quadric
sitting  in the 6-dimensional real ray space
$P^+(V)$ (diffeomorphic to $S^6$).  We can describe  points 
of the ray space  $P^+ (V)$ using homogeneous coordinates 
$[x] = [v,h] = [\lambda v, \lambda h]$,
$\lambda \in \R^+$ with $(v,h)\in Im\H\oplus\H =V$.
 $C$ is
defined  by the homogeneous equation $\|v\|^2=\|h\|^2.$ 
Using the $\R^+$ action, we   normalize  $\|v \| = 1$,
proving that  $C$ is diffeomorphic to
the  product of spheres  $S^2 \times S^3 \cong S^3 \times S^2 =\tilde Q$
which appears  in theorem 1.

  Given a point $\R^+ x = [x] \in C,$  set $$x^\perp=\{y\in V|\langle x,y\rangle=0\}, \quad
x^0=\{y\in V| xy=0\}.$$ Then
\begin{proposition}
 $$\R x\subset x^0\subset (x^0)^\perp\subset x^\perp\subset V,$$
and the dimensions are $1,3,4,6,7.$
\end{proposition}

\pf Use  the   definitions of  the split octonion product (eq (\ref{product})) and inner product.\qed

When we projectivize,  $x^\perp$ maps to the tangent plane $T_{[x]} C$ to $C$ at $[x]$, and $x^0$
maps to a 2-dimensional subspace $D_{[x]}\subset T_{[x]} C$. 
Letting $[x]$ vary over $C$ we have
defined a   rank 2 distribution $D\subset TC$. 
This construction of $(C, D)$   depends only on
 on the algebraic structure of   $\tilde \O$, so that  $G_2=Aut(\tilde\O)$
 acts on $C$ preserving $D$.

\begin{proposition}
\label{projective}
The (ray) projective quadric $C$ is a homogeneous
space for $G_2$.  $C$  is diffeomorphic to $\tilde Q =S^3 \times S^2$
of theorem 1,  and is naturally
endowed with a   $G_2$-invariant  distribution  $D$ of rank $2$.
Viewed as a $G_2$-homogeneous space, the data for $(C,D)$ 
 coincides  with the data  $(G_2, P, W)$ of  section 3.3.   Viewed
 as a $K$-homogeneous space, its Lie algebraic data
 coincides with that of the rolling distribution $(\k, \h, D)$ for 
 the ratio $1:3$.   The distribution on $C$
pushes down to the rolling distribution for ratios $3:1$
under the two-to-one cover   $C = S^3 \to S^2 \to Q = SO_3 \times S^2$
\end{proposition} 

This proposition immediately implies
that part of the theorem we are going to prove:
that $G_2 \subset Aut(\tilde Q, \tilde D)$. 

{\bf Steps of the proof.}
In the paragraph   preceding the   
 proposition we proved that $C$ is a homogeneous
 space for  $G_2$, that $D$ is invariant under this $G_2$ action,
 and that $C$ is diffeomorphic to $\tilde Q$.
 Next, we will prove that the  coincidence of the $\g_2$-data for $(C,D)$ 
 and the data   $(\p, W)$ of the previous section.
 For this we will  use  the weights for the 
 $G_2$-representation space $V = Im(\tilde \O)$.

{\bf   Weights for the $7$-dimensional representation.}
 
 Here is  the weight diagram  for this representation. 
 
\begin{center}
\unitlength=1cm
\begin{picture}(5,4)(-2,-2)

\put(0,0){\circle{.3}}

\put(0,0){\circle{.18}}

\put(-.5,-.866){\circle{.25}}
\put(-.8,-1){X}
\put(.5,.866){\circle{.25}}
\put(-.5,.866){\circle{.25}}
\put(-.68,.866){X}
\put(.5,-.866){\circle{.25}}

\put(-.5,-.866){\circle{.1}}
\put(.5,-.866){\circle{.1}}

\put(.5,.866){\circle{.1}}
\put(-.5,.866){\circle{.1}}

\put(-1,0){\circle*{.25}}
\put(1,0){\circle{.25}}
\put(-1,0){\circle{.1}}
\put(1,0){\circle{.1}}


\put(-1.4,-0.1){X}

\put(1.4,-.866){\circle{.25}}
\put(-1.5,.866){\circle{.25}}
\put(-1.9,  .8){ X  }

\put(1.4,.866){\circle{.25}}
\put(-1.5,-.866){\circle{.25}}
\put(-1.9,  -1){X}

\put(0,1.5){\circle{.25}}
\put(-.2,1.6){X}
\put(0,-1.5){\circle{.25}}
\put(-0.3,-1.7){X}

\put(0,-1.4){\line(0,1){3}}
\put(-1.5,  -0.866){\line(5,3){2.77}}
\put(-1.5,  0.866){\line(5,-3){2.77}}

\put(  -1, 0){\line(1  , 0   ){1.85}}
\put(-.5,-.866){\line(3,5){.96}}
\put(-.5,.866){\line(3,-5){.96}}

\end{picture}

Figure 4: Weights and roots  associated with the representation $V$.

\end{center}

The  weights of the representation $V$  form a subset of the roots of $\g_2$.
In figure 4 we redrew the root diagram of $\g_2$, marking  
those roots  which are  weights for   $V$   with bullseye's. 
They   are the six short roots and one zero root.
The corresponding weight spaces $V_w$ are all one-dimensional.
The black dot is a selected weight vector and corresponds to 
a `choice of base point' for $C$.  The meaning of the X's will be given below.

\bn{\bf A reminder of the meaning of the weight  diagram.} 

We   recall the general case of a  representation $V$ of a
semi-simple Lie algebra $\g$
with Cartan subalgebra $\t$.  
A {\it weight} for the  representation $V$ of $\g$
is an element $w \in \t^*$ such that there is a nonzero
vector $v \in V$ with the property that
$\zeta \cdot v = w(\zeta) v$ for all $\zeta \in \t$.    The space of $v$'s for a given weight $w$
is called the {\it weight space} for $w$ and is denoted $V_w$.  
If, for given $w \in \t^*$ there is no such nonzero $v$
then we set $V_w = 0$.   For a finite-dimensional representation
the set of weights is finite.       We have
$$V = \bigoplus_{w\in\t^*} V_w.$$
The roots of $\g$ are the non-zero weights of the adjoint representation. 

If, as in our situation, the roots for $\t$ are real, then its `torus' $T = exp(\t)$  is noncompact and  acts on the weight spaces by scaling, as follows. 
If $\lambda = exp(\xi) \in T$,  with $\xi \in \t$,
then $\lambda e_w = exp(w(\xi))e_w$ for $w \in V_w$.
 
From $\zeta \cdot \xi  \cdot v =  \xi \cdot \zeta   \cdot v + [\zeta, \xi ] \cdot v$
it follows that if $v$ is in the weight space for $w$
and $\xi \in \g_{\alpha}$ is in the root space for $\alpha$ then $\xi v$
is in the weight space for $w + \alpha$ (which, as above, could be zero).
In other words: 
$\g_{\alpha} V_w \subset V_{w + \alpha}$. This inclusion is  half of the rule:
\begin{equation}
w \hbox{ a weight}, \alpha \hbox{ a root} \Rightarrow  \g_{\alpha} V_w =  V_{w + \alpha}
\label{rule}
\end{equation}
which is true for $V$.
It follows in particular that if $v \in V_w$
and $\xi \in \g_{\alpha}$ and if $w + \alpha$ is not a weight for the representation,
then $\xi (v) = 0$.  We will use this fact momentarily.

We now construct  the weight spaces and the action of
the torus  for our $G_2$-representation
$V = Im(\tilde \O)$. 
Let $n$ be an imaginary quaternion.
Then $(n, n)$ and $(n, -n)$ are both null vectors in $V$.
 Take as  basis for $V$:
\begin{equation}
e_1 = \frac{1}{2}(i, i), e_2 = \frac{1}{2}(j, j), e_3 = \frac{1}{2}(k,k) ; f_1 = \frac{1}{2}(i, -i), f_2 = \frac{1}{2}(j, -j), f_3 = \frac{1}{2}(k, -k)
\label{basis}
\end{equation}
and 
$$U = (0,1).$$  Then we have
the multiplication table:
$$e_i ^2 = f_i ^2 = 0 $$
$$e_i f_j = f_j e_i = 0,  \hskip .2cm \hbox{ if }  i \ne j$$
$$e_i e_j = f_k;  i, j, k \hbox{ a  cyclic permutation of } 1,2,3 $$
$$f_i f_j =   e_k  ;   i, j, k \hbox{ a  cyclic permutation of } 1,2,3 $$
 $$e_i f_i =   -\frac{1}{2} + \frac{1}{2}U$$
 $$f_i e_i =  -\frac{1}{2}  - \frac{1}{2}U$$
 $$e_i U =  e_i$$
 $$f_i U = -  f_i $$
 To complete the multiplication table, use that the conjugate
 of  $xy$ is  $\bar x \bar y$,
 so that if $x, y \in V = Im(\tilde \O)$ we have 
 $y x = \bar z$ where
 $z =  xy$.  
 Thus, for example since $\bar f_k = -f_k$
 we see that  $e_j e_i =-f_k$,  for   $i,j,k$ a cyclic permutation of $1,2,3$. 
Now let $\lambda_1, \lambda_2, \lambda_3$
be nonzero reals with $\lambda_1 \lambda_2 \lambda_3 = 1$.
Let $\alpha_i, \beta_i, \gamma_i$
and  $\tilde \alpha_i, \tilde \beta_i,  \tilde \gamma_i$ be real exponents for $i = 1, 2, 3$
satisfying $\alpha_i + \beta_i + \gamma_i = 0$
Then the scaling  transformation
$$e_i \mapsto \lambda_1 ^{\alpha_i} \lambda_2 ^{\beta_i} \lambda_3 ^{\gamma_i} e_i$$
$$f_i  \mapsto \lambda_1 ^{\tilde \alpha_i} \lambda_2 ^{\tilde \beta_i} \lambda_3 ^{\tilde\gamma_i}  f_i$$
together with $z \mapsto z$ 
preserves the multiplication table, and hence defines an element of   $G_2$,
provided 
$$\tilde \alpha_i = - \alpha_i,  \tilde \beta_i = - \beta_i, \tilde \gamma_i = - \gamma_i$$
and provided that 
$(\alpha_i, \beta_i, \gamma_i)$ are multiples of the values from the following
weight table 

$$\begin{array}{|c||c|c|cl}
\hline
 &\alpha_i&\beta_i&\gamma_i\\
\hline\hline
i = 1  &2 &-1& -1\\
\hline
i=2& -1 &2 &-1  \\
\hline
i=3  &-1 & -1& 2   \\
\hline
\end{array}
$$
These scaling transformations generate the Cartan $T$ of $G_2$,
and the table  gives the corresponding weights.  Thus for example
$e_1$ is a weight vector with corresponding weight, relative to
the basis for $\t$,  being $(2,-1,-1)$.  Here we view $\t$
as being the collection of vectors $(a, b, c)$ with $a + b + c = 0$.
  Looking at the inner products of these vectors we see that they are arranged
on the weight diagram according to:

\begin{center}
\unitlength=1cm
\begin{picture}(5,4)(-2,-2)

\put(0,0){\circle{.3}}

\put(0,0){\circle{.18}}

\put(-.5,-.866){\circle{.25}}
 \put(.5,.866){\circle{.25}}
\put(-.5,.866){\circle{.25}}
 
\put(.5,-.866){\circle{.25}}
\put(.8, 1){$e_2$}
\put(.76, -1){$e_3$}

\put(-.5,-.866){\circle{.1}}
\put(.5,-.866){\circle{.1}}
\put(-.8, -1.18){$f_2$}
\put(-.8, 1.1){$f_3$}

\put(.5,.866){\circle{.1}}
\put(-.5,.866){\circle{.1}}

\put(-1,0){\circle*{.25}}
\put(1,0){\circle{.25}}
\put(-1,0){\circle{.1}}
\put(1,0){\circle{.1}}

\put(-1.45, 0){$e_1$}
\put(1.2, 0){$f_1$}


 \put(.2, .2){$z$} 
 
\put(  -1, 0){\line(1  , 0   ){1.85}}
\put(-.5,-.866){\line(3,5){.96}}
\put(-.5,.866){\line(3,-5){.96}}

\end{picture}

Figure 5:  The weight space basis. 

\end{center} 

We are now in a position to compute  the $\g_2$-data
associated to $(C, D)$ from the proposition. 

{\bf Weight vectors are null vectors: }   Because the inner product is $G_2$-invariant,
the $\g_2$ action on $V$ satisfies $\langle \xi x,  x \rangle = 0$
for any $\xi \in \g_2$, $x \in V$. Take $x$ a weight vector with nonzero weight $w$,
and   
take  $\xi \in \t$ with $w(\xi) \ne 0$.  From 
$\langle \xi x,  x \rangle =  w(\xi) \langle x, x \rangle$
we have that   $x$ is a null-vector.

{\bf Computing the isotropy data.}
Set $c_0 = [e_1]$, the ray through $e_1$.
We must show that the isotropy group of $c_0$
is $P$.  

We begin by showing that  the isotropy algebra
$\g_{c_0}$  of $c_0 = [e_1] $ is  
$\p$.   We have that   $\g_{c_0} = \{ \xi \in \g_2 :  \xi e_1 = \lambda e_1 \hbox{ for some
real number } \lambda \}$.
 The black dot in figure 4 indicates
the weight space spanned by $e_1$,  
with corresponding   weight 
by $w_1$.  (This weight is the root
marked $- \sigma_3$ in Figure 3.) 
According to the addition rule, 
(\ref{rule}) if $\alpha \in \t^*$ is a root
and $w_1 + \alpha$  {\it is not a weight} for $V$,
then $\xi_{\alpha} x_0 = 0$.    Those
roots for $w_0 + \alpha$ {\it is not a weight} are marked by X's in figure 4.
The sum of the corresponding $\g_{\alpha} \subset \g_2$ is a vector
space of elements $\xi$ satisfying   $\xi (e_1) = 0$.
Now  the isotropy algebra $\g_{c_0}$of the ray through $e_1$
consists of all those $\xi$ such that $\xi e_1 = \lambda e_1$
for some real scalar $\lambda$.   
The elements  $ H \in \t$   act on $e_1$ by
scalar multiplication by $\lambda = w_1 (H)$.
Referring to the diagram then, we see that
$\p \subset \g_{c_0}$.   But there is no
subalgebra of $\g_2$ lying between $\p$ and all of $\g_2$.  
It follows that the isotropy algebra for the ray {\it is} $\p$.

  It follows  from this Lie algebra computation 
    that the isotropy subgroup $G_{c_0} $
{\it contains} $P$ and has Lie algebra  equalling the Lie
algebra $\p$ of $P$.  Now   $P$
was defined to be the {\it connected} Lie subgroup of $G_2$
whose Lie algebra is $\p$, thus to show
$G_{c_0} = P$ it suffices to show that $G_{c_0}$ is connected.
We demonstrate connectivity 
by applying  the homotopy exact sequence
to the fiber bundle $G_{c_0} \to G_2 \to C = G_2 / G_{[x]}$.
This exact sequence
is   $\ldots  \to \pi_1(C) \to \pi_0 ( G_{c_0}) = \pi_0 (G_2) \to \pi_0 (C)$.
Since $C$ is simply connected and connected we
get that  $\pi_0 ( G_{c_0} = \pi_0 (G_2)$
and since $\pi_0 (G_2) = 0$ we have our connectivity:
$\pi_0 (G_{c_0}) = 0$.  

We have   established the isotropy ($P$ part) of the data for $(C, D)$.

{\bf Computing the distribution data}.    
The distribution plane  $D(c_0)$   at $c_0$ 
corresponds to $e_1 ^0$ -- the 
subspace $S \subset V$ consisting  of
those  vectors $y \in V$ for which 
$e_1 y = 0$.    From the multiplication table 
following the description of  our basis (\ref{basis})  we see that 
$S = span \{ e_1, f_2, f_3 \}$.
From Figure 4, we see that weights
corresponding to $f_2, f_3$, , say $w_2, w_3$,
are given by $w_2 = w_1 + \sigma_1$,
$w_3 = w_1 + (-\sigma_2)$.  Compare Figure 3.  
Let $x_1, y_1 \in \g_2$ be the corresponding
nonzero root vectors for $\sigma_1, - \sigma_2$. 
(We follow the $x, y$ notation from Figure 5, Appendix 2.)
It follows from   rule (\ref{rule})
that  $f_2$ is a multiple of   $x_1 (c_0)$
and $f_2$ is a multiple of $y_2 (c_0)$.
In other words, $S = W(e_1)$ mod $\p (e_1)$
where $W$ is the space spanned by the roots indicated
by the pluses in Figure 2.  We have proved
that the Lie algebraic   data for
$(C, D)$ is  $(\g_2, \p, W)$.

{\bf The covering map.} 
On the Lie algebra level we have shown that the data for
$(C, D)$ is  $(\g_2, \p, W)$. As computed in section 3.3, Cor. 1, upon
restricting the action of $G_2$  to $K$ this Lie algebraic data
 $(\g_2, \p, W)$ corresponds to  
the data $(\so_3 \oplus \so_3, \h, D(3;1))$.
Thus, up to a finite cover, $(C,D)$ is the rolling distribution.
Now $C$ is simply connected, and  $2:1$ covers  $Q$.
This covering map  $C = S^3 \times S^2 \to Q = SO_3 \times S^2$
is realized by forming the quotient of
$C$ by the $\Z_2$ subgroup
generated by image of  $\sigma = (\pm 1, 1) \in K = SU_1 \X_{\pm(1,1)} SU_2$.
 Being an element of the symmetry group   $\sigma$ preserves
the distribution $D$ on $C$, and so $D$ does  push down to
the rolling space $Q$. The $\k$ data of the pushed-down distribution
 remains   $(\so_3 \oplus \so_3, \h, D(3;1))$. Thus the pushed-down distribution
is isomorphic to the rolling distribution on $Q$.  

QED

\section{Summary. Lack of  action on the rolling space. The theorem is done.    }

We have proved that $G_2 \subset Aut( \tilde Q, \tilde D)$
and that is all that we are going to prove of
theorem 1,  with the exception of   the fact that the 
  $G_2$-action  does {\it not} descend to  $Q$.  (Recall from the introduction
  we are not going to prove that $G_2 = Aut( \tilde Q, \tilde D)$.)
  To prove that the $G_2$ action does not descend to $Q$, we 
  realize as above that   $Q = \Z_2 \backslash C$ 
where  the  $\Z_2 \subset K \subset G_2$ is generated by $\sigma = (\pm 1, 1)$.
Now we use the following fact about group actions.  Suppose
that a group $G$ (here $G_2$) acts effectively on a set $C$
and that $\Gamma \subset G$.   (``Effectively'' means that the only
group element acting as the identity on $C$ is the identity.)  
Then the action of an element $g \in G$ descends to the quotient space $\Gamma \backslash C$ if and only if $g \Gamma g^{-1} = \Gamma$.  
In particular, if $\Gamma$ is not normal in $G$ then
the action of all of $G$ does not descend to the quotient $\Gamma \backslash C$. 
Returning to our situation, we see that  
if the $G_2$ action were to descend then this
$\Z_2$ generated by $\sigma$ 
would have to be normal.   But a discrete normal subgroup
of a connected Lie group is central, and $G_2$ has
no center.  See Appendix A, or \cite{Vogan}.  So our $\Z_2$ is not normal,
and the $G_2$ action does not descend.

{\bf Remark.}  Had we used  lines instead of rays
when constructing  $C = \tilde Q$, we would have
arrived at a quadric  $Q_f$ in the standard real   projective space $P(V)$
which is double covered by $C = \tilde Q$.
(The subscript `f' is for `false'.) 
$Q_f \subset P(V)$ is diffeomorphic to    $S^3 \times_{\Z_2} S^2 = {\pm I} \backslash C$ where
the notation $\times_{\Z_2}$
indicates that we divide out by the action of the involution $(v, h) \mapsto (-v, -h)$.
(This involution does not lie in $G_2$.)
$C = \tilde Q$ double-covers both  $Q_f$  and
$Q$, and the distribution $\tilde D$ pushes down to both covered spaces. 
 But $Q_f$ is   topologically distinct  from $Q$.   Both
 $Q$ and $Q_f$ are $\SO_3$-bundles over $S^2$. $Q$ is  the trivial $SO(3)$-bundle. 
$Q_f$ is the other one. (Since $\pi_1 (\SO_3) = \Z_2$ there are precisely
two topologically distinct $\SO_3$ bundles over $S^2$.)
Because $-I \in GL(V)$ commutes with the $G_2$ action on $V$
the $G_2$-action on $\tilde Q$ does descend to $Q_f$. 
We find it  curious that   the action of $G_2$
on $\tilde Q$ does descend to this `false'  rolling configuration space $Q_f$,  but not to the real one $Q$.

 \appendix{}
 \section{Covers. Two   $G_2$'s. }
\label{twoG2s}
 
To understand our results, it helps to understand that  
 up to isomorphism,  there are precisely two connected $G_2$'s:
 the adjoint one which is the one we have been using,
 and the simply connected one, which is the universal cover
 of the adjoint one.  
 For a general semi-simple Lie algebra $\g$  we can always form
the simply connected Lie group $\tilde G$ having $\g$ as its   Lie algebra.
If $Z$ is the center of  $\tilde G$, then $Ad(\tilde G) = \tilde G/Z$
where $Ad(\tilde G)$ is the image of $G$ under the adjoint map
from $\tilde G$ to $Hom(\g)$.  If $Z \ne I$ then $G \ne Ad(G)$.
There are as   many distinct connected Lie groups with
algebra $\g$ as there are distinct subgroups of $Z$, these being
the connected topological groups covered by $\tilde G$ and covering   $Ad(G)$. 
So,  when  $Z = \Z_2$ there are precisely two such Lie groups,   $\tilde G$, the simply connected one, and $G = Ad(\tilde G)$, the adjoint one.  

We find on p. 3 of Vogan \cite{Vogan} that the center of the simply connected
$G_2$ is  indeed $\Z_2$, and hence we have precisely two $G_2$'s.   
It will be useful to explain a few details of this computation of $Z(G_2)$.
The universal cover of any $G$ contracts onto its maximal compact.
Thus, if the maximal compact of $Ad(G)$ has finite fundamental group,
then the universal cover $\tilde G \to Ad(G) = \tilde G /Z$
is a finite cover, and so  the center $Z$ must be  finite.
(At the other extreme,  the maximal compact of 
$SL(2, \R)$ is a circle group, corresonding to the fact that its 
universal cover has infinite center  $\Z$.)   

We saw above  that the  Lie algebra of 
the maximal compact of  any   $G_2$ realizing $\g_2$
 is  $\k = \so_3 \times \so_3$.  The  connected Lie groups $K$ 
  having $\k$ as Lie algebra have fundamental groups 
 consisting of
  either $1$, $2$ or $4$ elements.  
  It follows that the center $Z(G_2)$ of any $G_2$   is finite, and 
hence compact.  Being compact and central, this center lies
in every maximal compact:  $Z(G_2) \subset K \subset G_2$. 
If we take the simply connected  $G_2$, call
it $\tilde G_2$,  then its maximal compact is 
$\tilde K = \SU_2 \times \SU_2$. The center of  $\tilde K$   is  the group of the four elements   $(\pm 1, \pm 1)$.  
The center of $\tilde K$ need not be the center of $G_2$
but it must 
contain it:   $Z(\tilde G_2) \subset Z(\tilde K)$.
To see what the actual center of $\tilde G_2$  is, it suffices to see how  
$Z(\tilde K) \subset \tilde K$ acts on 
the Lie algebra $\g_2$  under the adjoint action.
This can be done using roots. 
The center of $\tilde G_2$ is that part of $Z(\tilde K)$ which
acts trivially on $\g_2$. A computation using roots and
the restriction of the adjoint representation to $\tilde K$ shows that
this part is $(1,1)$ and $-(1,1)$.

\section{The isomorphism of $\k$ and $\so_3 \oplus \so_3$ from Proposition 3. }
\label{maxcpt}

We complete the  proposition 3 from section 5, in which the
explicit identification of  $\so_3 \oplus \so_3$
as the Lie algebra   $\k$ of the maximal compact in $\g_2$.   We  follow  Serre  \cite{Serre}, 
page VI-11: $\g_2$ is Lie-generated by the elements
$x,y,h,X,Y,H$, subject to the following relations, which one can read
off the root diagram.

$$\begin{array}{llll}
[x,y]=h,& [h,x]=2x,& [h,y]=-2y,& \\
\left[X,Y\right]=H,& [H,X]=2X,&[H,Y]=-2Y;&\\
\left[h,X\right]=-3X, &[h,Y]=3Y;&[H,x]=-x,&[H,y]=y;\\
\left[x,Y\right]=[X,y]=[h,H]=0;&&&\\
\left[ad(x)\right]^4 X=0; &\left[ad(X)\right]^2x=0;&&\\
\left[ad(y)\right]^4Y=0;& \left[ad(Y)\right]^2x=0.&&\\
\end{array}
$$

Taking Lie brackets of the vectors $x,y,h,X,Y,H$ we generate a complete set  $\{x_i,X_i, y_i,Y_i|i=1,2,3\}$ of root vectors for
$\g_2$, which, together with the basis $h,H$ for the Cartan subalgebra form a basis for $\g_2$ as follows:
$$\begin{array}{llllll}
x_3=x,&X_1=X, &x_2=[x,X_1],&x_1=[x,x_2],&X_2=[x,x_1],&X_3=[X_1,X_2];\\
y_3=y,& Y_1=Y,&y_2=-[y,Y_1],&y_1=-[y,y_2],&Y_2=-[y,y_1],&Y_3=-[Y_1,Y_2].
\end{array}
$$
We label each root in the diagram
with the corresponding root vector.

\bn

\begin{center}

\unitlength=1cm
\begin{picture}(5,4)(-2,-2)
\put(0,0){\circle*{.4}}
\put(.5,.866){\circle*{.25}}
\put(.65,1.05){$x_1$}
\put(-.5,-.866){\circle*{.25}}
\put(-.75,-1.25){$y_1$}
\put(-.5,.866){\circle*{.25}}
\put(-.95,1.05){$x_2$}
\put(.5,-.866){\circle*{.25}}
\put(.65,-1.25){$y_2$}
\put(1,0){\circle*{.25}}
\put(1.2,-0.1){$x_3=x$}
\put(-1,0){\circle*{.25}}
\put(-2.4,-0.1){$y=y_3$}
\put(0,1.7){\circle*{.25}}
\put(-.2,2){$X_3$}
\put(0,-1.732){\circle*{.25}}
\put(-.2,-2.3){$Y_3$}
\put(-1.5,.866){\circle*{.25}}
\put(-3.1, 1.05){$X=X_1$}
\put(1.4,-.866){\circle*{.25}}
\put(1.6,  -1.25){$Y_1=Y$}
\put(1.4,.866){\circle*{.25}}
\put(1.7, 1.05){$X_2$}
\put(-1.5,-.866){\circle*{.25}}
\put(-2.2,-1.25){$Y_2$}

\put(0,-1.6){\line(0,1){3.3}}
\put(-1.4,  -0.8){\line(5,3){2.77}}
\put(-1.5,  0.866){\line(5,-3){2.77}}
\put(  -.9, 0){\line(1  , 0   ){1.85}}
\put(-.46,-.75){\line(3,5){.96}}
\put(-.5,.866){\line(3,-5){.96}}

\end{picture}

\vskip .8cm
Figure 5.  A basis for the Lie algebra.
\end{center}

\bn

We end up with a ``nice'' basis wrt which the structure constants are
particulary pleasant; they are integers and have symmetry properties
which facilitate greatly the work involved in their determination; you
can also apply some elementary $\sl_2$ representation theory that
further facilitate the calculation; it helps to work with the root
diagram nearby.

\bn{\bf Symmetry properties of the structure constants.}  Suppose $\alpha,\beta$ are two roots such
that $\alpha+\beta$ is also a root. Let $E_\alpha, E_\beta$
be the corresponding root vectors, as chosen above. Then
$[E_\alpha,E_\beta]=c_{\alpha,\beta}E_{\alpha+\beta}$, for some
non-zero constant $c_{\alpha,\beta}\in\Z$. The nice feature of
our base is that the structure constants satisfy
$$c_{-\alpha,-\beta}=-c_{\alpha,\beta}.$$ This cuts in half the amount
of work involved, since you need only consider say $\alpha>0$ (the positive roots are the six dots in the last root diagram marked 
with $x$'s and $X$'s).
Combining this with the obvious $c_{\alpha,\beta}=-c_{\beta,\alpha}$
(antisymmetry of Lie bracket) you obtain
$$c_{\alpha,-\beta}=c_{\beta,-\alpha}.$$
This cuts in half again the amount of work.

\begin{proposition} The structure  constants of $\g_2$, with respect to the basis of root vectors $\{x_i,X_i,y_i, Y_i | i=1,2,3\}$
and the Cartan algebra elements $\{h,H\}$ are given as follows. The
basis elements are grouped in three sets: positive (three $x$'s and
three $X$'s), negative (three $y$'s and three $Y$'s), and Cartan
subalgebra elements ($h$ and $H$).
\begin{itemize}
\item $[Positive, positive]$: other then the ones given above, and
those which are zero for obvious reasons from the root diagram (sum of
roots which is not a root): $$[x_1, x_2]= X_3.$$

\item $[ Positive, negative ]$:

$$\begin{array}{|c||c|c|c|c|c|c|}
\hline
c_{\alpha,\beta}&y_1&y_2&y_3&Y_1&Y_2&Y_3\\
\hline\hline
x_1&1&4&-4&0&12&-12\\
\hline
x_2& 4 &1 &-3 &1 &0 &3 \\
\hline
x_3  &-4 & -3& 1& 0&-3&0  \\
\hline
X_1 &0 &1 &0 &1 &0&-1 \\
\hline
X_2 & 12& 0& -3& 0& 1&36\\
\hline
X_ 3& -12& 3& 0&-1 &36&1\\
\hline
\end{array}
$$

 \mn The 1's on the diagonal stand for the relations $[x_i,y_i]=h_i,$
 $[X_i,Y_i]=H_i,$ where, in terms of our basis $\{h,H\}$ for the
 Cartan subalgebra, $$h_1=8h+12H,\quad h_2=h+3H, \quad h_3=h,$$
$$ H_1=H,\quad H_2=36(h+H), \quad H_3=36(h+2H).$$
\item $[Cartan, anything]$: this is coded directly by  the root diagram:

- $ad(x)$  has eigenvalues and eigenvectors
$$\begin{array}{l||c|c|c|c|c|c|c}
\hbox{eigenvalue}&3&2&1&0&-1&-2&-3\\
\hline
\hbox{eigenvectors}&X_2,Y_1&x_3& x_1, y_2&X_3,Y_3, h,H&x_2,y_1&y_3&X_1,Y_2
\end{array}$$

- $ad(X)$ has
eigenvalues and eigenvectors
$$\begin{array}{l||c|c|c|c|c}
\hbox{eigenvalue}&2&1&0&-1&-2\\
\hline
\hbox{eigenvectors}&X_1&X_3,x_2,y_3,Y_2&x_1,y_1,h,H&X_2,x_3,y_2,Y_3&Y_1
\end{array}$$

\end{itemize}
\end{proposition}

\pf This is elementary, using only the Jacobi identity, but takes time. We will give as a typical example
the calculation of $[x_1, x_2]$:

$$\begin{array}{rcll}
[x_1,x_2]&=&[x_1,[x,X]]&\quad\hbox{(by definition of $x_2$)}\\
&=&[x,[x_1,X]]+[X, [x,x_1]]&\quad\hbox{(Jacobi identity)}\\
&=&[X, [x,x_1]]&\quad\hbox{(since $[x_1,X]=0$)}\\
&=&[X, X_2]=X_3&\quad\hbox{(by definitions of $X_2, X_3$).}
\end{array}
$$
The rest of the relations are derived in a similar fashion. \qed

\bn

Now we are ready to define the generators of the Lie algebra of a maximal compact subgroup $K\subset G_2$. Let
$$L_1=X_1-Y_1, \quad L_2={X_2-Y_2\over 6}, \quad  L_3={X_3-Y_3\over 6},$$
$$S_1= {x_1-y_1\over 4}, \quad S_2= {x_2-y_2\over 2}, \quad S_3= {x_3-y_3\over 2}.$$

Using the commutation relations of the last Proposition  one checks easily  that
$$[L_i,L_j]=\epsilon_{ijk}L_k,\quad [L_i,S_j]=\epsilon_{ijk}S_k,\quad [S_i,S_j]=\epsilon_{ijk}({3\over 4}L_k-S_k).$$

\bn  Note: the strange-looking   coefficients 2,4,6 in the definition of the $L_i, S_i$ are chosen precisely so that we get
these pleasing commutation relations.

\section{The rolling distribution in Cartan's thesis}

\subsection{Cartan's constructions and claims.} 
\label{claims}

In E. Cartan's thesis \cite{C2}, p.146,  we find the following constructions:
consider $V=\R^7=\R^3\X \R^3\X \R$ with coordinates $(\x,\y,z)$, where $\x,\y\in\R^3$, $z\in\R$,  
and the following 15 linear vector fields (hence linear operators) on  $V $:

\begin{itemize}
\item $X_{ii}=-x_i \partial_{x_i} + y_i\partial_{y_i}+
{1\over 3}\sum_{j=1}^3 (x_j\partial_{x_j}-y_j\partial_{y_j})$, $i=1,2,3$.

\item $X_{i0}=2z\partial_{x_i}- y_i \partial_z - x_j \partial_{y_k} + x_k\partial_{y_j},$ $(ijk)\in A_3=\{(123),(231),(312)\}.$

\item $X_{0i}=-2z\partial_{y_i}+ x_i \partial_z +y_j \partial_{x_k} - y_k\partial_{x_j},$ $(ijk)\in A_3.$

\item $X_{ij}=- x_j \partial_{x_i} + y_i\partial_{y_j},$ $i\neq j,$ $i,j=1,2,3.$

\end{itemize}

Cartan makes the following claims without proof:  

\begin{enumerate}

\item The linear span of these 15 operators is a 14 dimensional Lie subalgebra $\g\subset\End(V)$ isomorphic to $\g_2$. 

\item $\g$ preserves the  quadratic form on $V$ given by 
$$J=z^2+\x\cdot\y.$$

\item The linear group $G\subset\GL(V)$ generated by $\g$ acts transitively on the projectivized null cone of $J$.

\item  $G$ preserves the system of 6 Pfaffian  equations  on $V,$
given by the 6 components of 

$$
\left\{
\begin{array}{rcl}
\alpha&:=&zd\x-\x dz+\y\X d\y=0,\\
 \beta &:=&zd\y-\y dz+\x\X d\x=0,
 \end{array}
 \right.
 $$
which have as a consequence 
$$
\left\{
\begin{array}{rcl}
\gamma_1&:=&zdz+\x\cdot d\y=0,\\
\gamma_2&:=&zdz+\y\cdot d\x=0.
 \end{array}
 \right.
 $$
 
 $$$$
\item $G$ preserves  a 5 parameter family of 
3 dimensional linear subspaces of $V,$ contained in the null cone of $J$, 

$$
\left\{
\begin{array}{l}
\x-z\a +\b\X \y=0,\\
\y-z\b +\a\X \x=0,
 \end{array}
 \right., 
 $$
where 
$$\a\cdot\b+1=0.$$

\end{enumerate}

Our goal in this appendix is to sketch proofs of these claims,
provide a minor correction in one place,  
relate Cartan's  construction to the  octonions,
and show how they contain, in essence, the 
construction of the rolling distribution $\tilde Q$
via projective geometry, as in the proposition \ref{projective}
from section 5.

\subsection{Relation with Octonions}.
\label{connectOctonions}

Recall the basis $e_i ,   f_i , U $ of section \ref{basis} for $V$ (imaginary split octonions) 
with its consequent multiplication table. 
Make the change of basis $e_i \mapsto -e_i$,
keeping $f_i, U$ as they were,  thus changing the
signs of some entries of the multiplication table.
Use this new   basis   $E_i = -e_i ,   f_i , U$
  to identify   $V$ with $\R^3 \times \R^3 \times \R$
by setting   $(\x, \y, z)  = \Sigma x_i E_i + \Sigma y_i f_i + z U \in V$. 
 Referring to the multiplication
table we compute 
\begin{eqnarray*}
(\x, \y, z) (\x', \y', z') &= &(-\y \times \y' -z \x' + z' \x, \x \times \x' + z \y' - z' \y, \frac{1}{2}( \x \cdot \y' - \x' \cdot \y)) 
\\
&&+ 1 \{ z z' + \frac{1}{2} ( \x \cdot \y' - \x' \cdot \y) \}.
\end{eqnarray*}
The last term is in the real part of the split octonions, and not in $V$. 
It follows from this formula that $(\x, \y, z)^2 =J$, of Cartan's claim
2 in the preceding paragraph.
Multiplying out   $(\x, \y, z) (d \x, d \y, dz)$ we  find that
$$(\x, \y, z) (d \x, d \y, dz)
 = (\alpha, \beta, \frac{1}{2} (\gamma_1 - \gamma_2)) + 
 1 \{  \frac{1}{2} (\gamma_1 +  \gamma_2) \},$$
 where  $\alpha, \beta, \gamma_1, \gamma_2$
 are as in   Cartan's claim 4 of the previous paragraph. 
 It follows that the any element of $G_2 = Aut(\tilde \O)$
 preserves $J$ and preserves   the Pfaffian system of Cartan's claim 4.  
The distribution  $D$ defined by this system  is, upon restriction
to the   null cone $\{J =0 \} \setminus \{ 0 \}$,  precisely the    distribution $D$  which  we defined  in the final section of
our paper:  $D(\x, \y, z) : = 
\{ (\a, \b , c) :  (\x, \y, z) (\a, \b, c) = 0 \}$.  
It follows that Cartan's construction, 
pushed down to the space of rays using the
$\R^+$-action,   yields
precisely our $\tilde Q$.  

\subsection{Commentary  and proofs of Cartan's claims.}

\subsubsection{Definition of $\g_2$} 

{\bf The Cartan subalgebra.} The first 3 operators of claim 1  are linearly dependent since
$\sum_i X_{ii}=0$. This is the only linear relation (proof below) and explains why $\g$ is 14 dimensional and not 15 dimensional.  The flows
of $3 X_{ii}$ generate the scalings
$x_i \mapsto \lambda_1 ^{\alpha_i} \lambda_2 ^{\beta_i} \lambda_3 ^{\gamma_i}x_i$
$y_i \mapsto \lambda_1 ^{-\alpha_i} \lambda_2 ^{-\beta_i} \lambda_3 ^{-\gamma_i}y_i$
$z \mapsto z$ as described in section 5.  Hence these operators  should span the
Cartan $\t$ of $\g = \g_2$.

\begin{proposition} $\g$ is a  14 dimensional   Lie subalgebra of $\End(V)$,  isomorphic to $\g_2$, with Cartan subalgebra as just described. 
\end{proposition}

\pf It is convenient to put $\g$ in block matrix form. For each ${\bf u}\in\R^3$ let  $R_{\bf u}\in\End(\R^3)$ be given by $ {\bf v}\mapsto{\bf u}\X {\bf v}$; i.e. 
$$R_{\bf u}=\left(\begin{array}{ccc}
0&-u_3&u_2\\
u_3&0&-u_1\\
-u_2&u_1&0
\end{array}\right).$$

\mn Define the linear  map $\rho:\sl_3(\R)\X\R^3\X\R^3\to\End(V)$ by

$$
\rho(A,\b,\c)=\left(\begin{array}{ccc}
A&R_{\c} &2\b\\
-R_\b&-A^t&-2\c\\
\c^t&-\b^t&0
\end{array}\right).
$$
Now $\rho$ is clearly   injective, hence its image is a 14 dimensional linear subspace of $\End(V)$. 
Denote the components of $A,\b,\c$ by $a_{ij},b_i,c_i$ (resp.), then  it is easy to check that 
$$\rho(A,\b,\c)=-\sum_{i,j}a_{ij}X_{ij}+\sum_i b_iX_{i0}+\sum_i c_iX_{0i}.$$
This shows that $\g$ is the image of $\rho$ and hence a 14 dimensional subspace of $\End(V)$. 

\bn To show that $\g$ is a lie algebra one calculates that $$[\rho(A,\b,\c),\rho(A',\b',\c')]=\rho(A'',\b'',\c''),$$
where
\begin{eqnarray*}
A''&=&[A,A']+3(\b\c'^t-\b'\c^t)- [\b\cdot \c'-\b'\cdot\c]I,\\
\b''&=&A\b'-A'\b-2\c\X \c',\\
 \c''&=&-A^t\c'+A'^t\c+2\b\X \b'.
\end{eqnarray*}\\

These formulae show that $\{\rho(A,0,0)|A\in\sl_3(\R)\}$ forms  a lie subalgebra of $\g$ isomorphic to $\sl_3(\R)$.  This
subalgebra corresponds to the sum of the long root
spaces in the root diagram, and  the Cartan subalgebra
(the sum of the $X_ii$) as identified earliter.
The formulae also show that 
the images of the   $\rho(0,\b,0)$  and $\rho(0,0,\c)$ are stable  under the adjoint action of the Cartan, hence they must correspond to the remaining short roots.   

A tedious computation now yields the root diagram and the structure
constants of $\g_2$.   
\bn 

\subsubsection{Invariance of $J$}
Let $G_2\subset\GL_7(\R)$ be the subgroup generated by $\g$.

\begin{proposition}  $J$ is $G_2$-invariant.  
\end{proposition}

\pf This is equivalent to showing that every $X\in\g$ is $J$-antisymmetric, i.e. that $X$ anti-commutes with 
$$
\left(\begin{array}{ccc}
0&I/2&0\\
I/2&0&0\\
0&0&1
\end{array}\right).
$$

One now checks easily that the set of $J$-antisymmetric matrices consists of the matrices of the form
$$\left(\begin{array}{ccc}
A&R_\c&2\tilde\b\\
-R_\b&-A^t&-2\tilde\c\\
\tilde\c^t&-\tilde\b^t&0
\end{array}\right),
$$
where $A\in\End(\R^3)$ and $\b,\tilde\b,\c,\tilde\c\in\R^3$. Looking at the formula for 
$\rho(A,\b,\c)$ we see that $\g$ is the subset of the $J$-antisymmetric matrices 
satisfying $\rm tr A=0, \b=\tilde\b,\c=\tilde\c$ (a codimension  7 condition).\qed

\subsubsection{Invariance of the Pfaffian system.}
 
\mn

\n {\bf Generalities.} A ``Pfaffian system''on  a manifold $M$ is given locally by the common kernels of a finite set of 1-forms, 
$$\alpha_1=\ldots=\alpha_m=0.$$
 Two sets of 1-forms $$\{\alpha_1,\ldots,\alpha_m\},\quad \{\beta_1,\ldots,\beta_n\},$$
give equivalent systems if one can express each element of one set as a linear combination 
(with coefficients in $C^{\infty}(M)$) of the elements of the other set. We write this as 
$$\alpha_i\equiv 0 \quad  \mod \beta_1,\ldots,\beta_n, \quad i=1,\dots, m,$$
and similarly for the $\beta$'s.

Consequently, if we want to prove that a system is preserved by some diffeomorphism $f:M\to M$
we must show that $$f^*\alpha_i\equiv 0\quad  \mod \alpha_1,\ldots,\alpha_m, \quad i=1,\dots, m,$$
and if we want to show that the flow of some vector field $X$ on $M$ preserves the system 
we must show that 
$$\L_X\alpha_i\equiv 0\quad  \mod \alpha_1,\ldots,\alpha_m, \quad i=1,\dots, m.$$

Given such a system we can consider the common kernels $D_x\subset T_xM$ of the 1-forms at each point $x\in M$. 
This is well defined independently of the 1-forms chosen to represent the system. If $\dim D_x$ (the rank of the system) 
is constant we obtain a distribution $D\subset TM$ (a subbundle of the tangent bundle). But the rank may vary. 
For example, the system on $\R$ given by $xdx=0$ has rank 1 at $x=0$ and rank 0 for $x\neq 0$.  However, if $G$ acts on $M$ preserving a Pfaffian system, then the rank must clearly be constant along the $G$-orbits.

\bn {\bf Cartan's Pfaffian system.    Rank jumps.  A  correction.}  Due to jumping of rank,
as discussed in the last remark, the  Pfaffian system
which Cartan defined by 
  the vanishing of  the 6 components of $\alpha,\beta$
 cannot be $G_2$ invariant, even when restricted to $\tC$,  the $J$ null cone.
For  at $(\e_1,0,0)$ the system   reduces 
to $dx_2=dx_3=dz=0$ and so has rank $4$. On the other hand, at
 the point $(\e_1,\e_2,0)$ the system is equivalent to  $dy_1=dx_2=dz-dy_3=dz+dx_3=0,$ and  so has rank  $3$. 
 And both points lie in  $\tC\setminus\{0\}$,
 which is a single $G_2$-orbit, contradicting $G_2$ invariance.  
A related  problem with Cartan's claim  4 of subsection (\ref{claims}) is his claim that $\gamma_1=\gamma_2=0$ is a consequence of $\alpha=\beta$. But this is true only on the $z\neq 0$ part of $\tC$. 

Both errors are fixed by  imposing  the extra equation $\gamma:=\gamma_1-\gamma_2=0.$  Then, as in section (\ref{connectOctonions}),
 we do obtain a $G_2$-invariant system on $V$.    Furthermore, 
 as proved immediately below, the two equations $\gamma_1=\gamma_2=0$ are indeed a consequence of $\alpha=\beta=0, \gamma=0$ on $\tC$, and are a consequence $\alpha=\beta=0$ on the subset $z\neq 0$ of $\tC$. 
 So   Cartan's claim is correct on the open dense set $z\neq 0$ of the null cone $\tC\subset V$. 
(See also page 11 of Bryant's paper on Geometric Duality \cite{Bryant}, where he adds the equation $\gamma=0$ to 
$\alpha=\beta=0$.)

\begin{proposition} The Pfaffian system on $V$ given by $\alpha=\beta=0, \gamma=0$ is  $G_2$-invariant. 
On $\tC$ the system is equivalent to $\alpha=\beta=0, \gamma_1=\gamma_2=0$. On the subset $z\neq 0$ of $\tC$
it is equivalent to $\alpha=\beta=0.$
\end{proposition}

\pf We prove the   claims of the last two sentences  first.
Note that  $\gamma_1+\gamma_2=dJ$. It follows that on  $\tC$, where $J=0$, 
we have that $\gamma_1=\gamma_2=0$ 
is a consequence of $\gamma:=\gamma_1-\gamma_2=0.$ Thus, restricted to $\tC$, the system $\alpha=\beta=0, \gamma=0$
is equivalent to $\alpha=\beta=0,\gamma_1=\gamma_2=0$. 
 Next, note that $\x\cdot\beta-\y\cdot\alpha=z\gamma.$ It follows that on $z\neq 0 $ the equation $\gamma=0$ is a consequence of  $\alpha=\beta=0$.\qed

It remains to establish invariance. 
We need to show that $$\L_X\alpha_i\equiv\L_X\beta_j \equiv\L_X\gamma\equiv 0 \mod \alpha_i,\beta_j,\gamma,$$
for all $X=\rho(A,\b,\c)\in\g.$ 
Divide into  3 cases, corresponding to $(A, 0,0)$ , $(0, \a, 0)$ and $(0, 0, \b)$
in our coordinatization of $\g$. 
\begin{itemize}
\item case 1: $X=\rho(A,0,0),$ $A\in\sl_3(\R)$.

\begin{lemma} If   $A\in\End(\R^3)$ and ${\bf u},{\bf v}\in\R^3$, then 
$$A({\bf u}\X {\bf v})+A^t{\bf u}\X {\bf v}+ {\bf u}\X A^t{\bf v}={\rm tr }A({\bf u}\X {\bf v}).$$
\end{lemma}
\pf Sketch: divide in 2 cases.  If  $A^t=-A$ then ${\rm tr }A=0$ and the identity is a consequence of the fact the 
$\SO_3$ preserves de cross product and that $\so_3$ are the antisymmetric matrices. If $A^t=A$ then can assume w.l.o.g. that $A$ is diagonal and do  an explicit easy calculation. \qed

Now since $$X(\x,\y,z)=(A\x,-A^t\y,0), \quad \alpha=zd\x-\x dz+\y\X d\y,$$
we get, using the lemma and  ${\rm tr }A=0$, that
\begin{eqnarray*}
\L_X\alpha&=&zAd\x-A\x dz-A^t\y\X d\y-\y\X A^td\y=\\
&=&A(zd\x-\x dz+\y\X d\y)=A\alpha\equiv 0\quad \mod \alpha.
\end{eqnarray*}
Similarly, $\L_X\beta=-A^t\beta\equiv 0\quad (\mod \beta).$

Finally, $\L_X\gamma=(A\x)\cdot d\y- \x\cdot (A^t d\y)=0.$

\item case 2: $X=\rho(0,\b,0),$ $\b\in\R^3$.

Here $$X(\x,\y,z)=(2\b z, -\b\X  \x, -\b\cdot\y),   $$
and one calculates that 
$$\L_X\alpha=\b\gamma,\quad \L_X\beta=\b\X\beta,\quad\L_X\gamma=\b\cdot\beta.$$

\item case 3: $X=\rho(0,0,\c), \c\in\R^3.$ The proof for this  case is very similar to the previous case.  Just interchange $\x$ and $\y$, and $\b$ and $\c$.
\end{itemize}
This completes the proof of invariance,
and hence the proof of the proposition.


\end{document}